\documentclass{siamart190516}

\usepackage{amssymb}
\usepackage{booktabs}
\usepackage{multirow}

\title{Jet Marching Methods for Solving the Eikonal Equation} \author{
  Samuel F. Potter \and Maria K. Cameron}

\renewcommand{\phi}{\varphi}

\newcommand{\m}[1]{\boldsymbol{#1}}
\newcommand{\nb}{\texttt{nb}}
\newcommand{\state}{\texttt{state}}
\newcommand{\valid}{\texttt{valid}}
\newcommand{\trial}{\texttt{trial}}
\newcommand{\far}{\texttt{far}}
\newcommand{\xhat}{\hat{\m{x}}}
\newcommand{\xlam}{\m{x}_{\m{\lambda}}}
\newcommand{\that}{\hat{\m{t}}}
\newcommand{\tlam}{\m{t}_{\m{\lambda}}}
\newcommand{\half}{1/2}

\newcommand{\lam}{\m{\lambda}}
\newcommand{\dphi}{\delta\hspace{-0.1em}\m{\phi}}
\newcommand{\mphi}{\m{\phi}}
\newcommand{\mpsi}{\m{\psi}}
\newcommand{\mell}{\m{\ell}}

\DeclareMathOperator{\Arg}{arg}
\DeclareMathOperator{\conv}{conv}

\DeclareMathOperator{\range}{range}

\begin{document}

\maketitle

\begin{abstract}
  We develop a family of compact high-order semi-Lagrangian
  label-setting methods for solving the eikonal equation. These
  solvers march the total 1-jet of the eikonal, and use Hermite
  interpolation to approximate the eikonal and parametrize
  characteristics locally for each semi-Lagrangian update. We describe
  solvers on unstructured meshes in any dimension, and conduct
  numerical experiments on regular grids in two dimensions. Our
  results show that these solvers yield at least second-order
  convergence, and, in special cases such as a linear speed of sound,
  third-order of convergence for both the eikonal and its gradient.
  We additionally show how to march the second partials of the eikonal
  using cell-based interpolants. Second derivative information
  computed this way is frequently second-order accurate, suitable for
  locally solving the transport equation. This provides a means of
  marching the prefactor coming from the WKB approximation of the
  Helmholtz equation. These solvers are designed specifically for
  computing a high-frequency approximation of the Helmholtz equation
  in a complicated environment with a slowly varying speed of sound,
  and, to the best of our knowledge, are the first solvers with these
  properties. We provide a link to a package online providing the
  solvers, and from which the results of this paper can be reproduced
  easily.
\end{abstract}



\section{Introduction}

Our goal is to develop a family of high-order semi-Lagrangian eikonal
solvers which use compact stencils. This is motivated by problems in
high-frequency room acoustics, although the eikonal equation arises in
a tremendous variety of modeling problems~\cite{Sethian:1999ab}.

In multimedia, virtual reality, and video games, precomputing room
impulse responses (RIRs) or transfer functions (RTFs) enables
convincing spatialized audio, in combination with binaural or surround
sound formats. Such an approach, usually referred to as
\emph{numerical acoustics}, involves computing pairs of RIRs by
placing probes at different locations in a voxelized domain,
numerically solving the acoustic wave equation, and capturing salient
perceptual parameters throughout the domain using a streaming
encoder~\cite{Raghuvanshi:2014aa,Raghuvanshi:2018aa}. These parameters
are later decoded using signal processing techniques in real time as
the listener moves throughout the virtual environment. Assuming that
the encoded parameters can comfortably fit into memory, a drawback of
this approach is that the complexity of the simulation depends
intrinsically on the highest frequency simulated. In practice,
simulations top out at around 1 kHz. The hearing range of humans is
roughly 20 Hz to 20 kHz, which requires these methods to either
implicitly or explicitly extrapolate the bandlimited transfer
functions to the full audible spectrum.

An established alternative to this approach is \emph{geometric
  acoustics}, where methods based on raytracing are
used~\cite{Savioja:2015aa}. Contrary to methods familiar from computer
graphics, the focus of geometric acoustics is different. Acoustic
waves are mechanical and have macroscopic wavelengths. This means that
subsurface scattering, typically modeled using BRDFs in raytracing for
computer graphics~\cite{Nicodemus:1965aa}, is less relevant, and is
limited to modeling macroscopic scattering from small geometric
features, since reflections from flat surfaces are specular in
nature. What's more, accurately modeling diffraction effects is
crucial~\cite{Schissler:2014aa}: e.g., we can hear a sound source
occluded by an obstacle, but we can't see it. A variety of other
geometric-acoustic methods exist beyond raytracing. Examples include
the image source method~\cite{Allen:1979aa} and frustum
tracing~\cite{Chandak:2008aa}.

Geometric acoustics and optics both assume a solution to the wave
equation based on an asymptotic high-frequency (WKB) approximation to
the Helmholtz equation~\cite{Popov:2002aa}. In this approximation, the
eikonal plays the role of a spatially varying phase function, whose
level sets describe propagating wavefronts. The prefactor of this
approximation describes the amplitude of these wavefronts. The WKB
approximation assumes a ray of ``infinite frequency'', suitable for
optics, since the effects of diffraction are limited. A variety of
mechanisms for augmenting this approximation with
frequency-dependent diffraction effects have been proposed, the most
successful of which is Keller's \emph{geometric theory of
  diffraction}~\cite{Keller:1962aa} (including the later \emph{uniform
  theory of diffraction}~\cite{Kouyoumjian:1974aa}).

The complete geometric acoustic field of multiply reflected and
diffracted rays can be parametrized by repeatedly solving the eikonal
equation, using boundary conditions derived from the WKB approximation
to patch together successive fields. A related approach is Benamou's
\emph{big raytracing} (BRT)~\cite{Benamou:1996aa,Benamou:1997aa}. This
approach requires one to be able to accurately solve the transport
equation describing the amplitude, e.g.\ using paraxial
raytracing~\cite{Popov:2002aa}. In order to do this, the first and
second order partial derivatives of the eikonal must be
computed. High-order accurate iterative schemes for solving the
eikonal equation exist~\cite{Zhang:2006aa,Xiong:2010aa,Luo:2014aa},
but their performance deteriorates in the presence of complicated
obstacles. Direct solvers for the eikonal equation allow one to
locally parametrize the characteristics (rays) of the eikonal
equation, which puts one in a position to simultaneously march the
amplitude. This enables work-efficient algorithms, critical if a large
number of eikonal problems must be solved.

Benamou's line of research related to BRT seems to have stalled due to
difficulties faced with caustics~\cite{Benamou:2003aa}. This is
reasonable considering that the intended use was seismic modeling,
where the eikonal equation is used to model first arrival times of
$P$-waves. In this case, the speed of sound is extremely complicated,
resulting in a large number of caustics~\cite{Versteeg:1994aa}. On the
other hand, in room acoustics, the speed of sound varies slowly. The
main challenge is geometric: the domain is potentially filled with
obstacles. This provides another motivation for compact stencils: such
stencils can be adapted for use with unstructured meshes, and the sort
of complicated boundary conditions that arise when using finite
differences are avoided entirely. In this work, our goal is to develop
the underlying approach to obtaining a compact higher-order
semi-Lagrangian eikonal solver. In future work, this will be applied
in the unstructured setting.

The solvers developed in this work are high-order, have optimally
local/compact stencils, and are label-setting methods (much like
Sethian's \emph{fast marching method}~\cite{Sethian:1999ac} or
Tsitsiklis's semi-Lagrangian algorithm for solving the eikonal
equation~\cite{Tsitsiklis:1995aa}). Additionally, being
semi-Lagrangian, they locally parametrize characteristics (acoustic
rays), making them suitable for use with paraxial
raytracing~\cite{Popov:2002aa}, the method of choice for locally
computing the amplitude. To the best of our knowledge, these are the
first eikonal solvers with this collection of properties.

We refer to our solvers as \emph{jet marching methods} to reflect the
fact that the key idea is marching the \emph{jet} of the eikonal (the
eikonal and its partial derivatives up to a particular
order~\cite{Shilov:1971aa}) in a principled fashion. Sethian and
Vladimirsky developed a fast marching method that additionally marched
the gradient of the eikonal in a short note, but did not prove
convergence results or provide detailed numerical
experiments~\cite{Sethian:2000aa}. Benamou and collaborators built on
these ideas by exploiting information about the eikonal equation to
obtain a compact upwind second-order finite difference scheme for
solving the eikonal equation~\cite{Benamou:2010vd}. Related methods
exist in the level set method community and are referred to as
\emph{gradient-augmented level set methods} or \emph{jet
  schemes}~\cite{Nave:2010aa,Seibold:2011aa}.

In the rest of this work we lay out these methods, providing detailed
numerical experiments. Our presentation is for unstructured grids in
$n$-dimensions, while our numerical experiments were carried out in
2D. We plan to extend these solvers to structured and unstructured
meshes in 3D and will report on these later in the context of room
acoustics applications.

\subsection{Problem setup}

Let $\Omega \subseteq \mathbb{R}^n$ be a domain, let $\partial\Omega$
be its boundary, and let $\Gamma \subseteq \Omega$. The \emph{eikonal
  equation} is a nonlinear first-order hyperbolic partial differential
equation given by:
\begin{equation}\label{eq:eikonal-equation}
  \begin{split}
    \|\nabla\tau(\m{x})\| &= s(\m{x}), \qquad \m{x} \in \Omega, \\
    \tau(\m{x}) &= g(\m{x}), \qquad \m{x} \in \Gamma.
  \end{split}
\end{equation}
Here, $\tau : \Omega \to \mathbb{R}$ is the \emph{eikonal}, a spatial
phase function that encodes the first arrival time of a wavefront
propagating with pointwise \emph{slowness} specified by
$s : \Omega \to (0, \infty)$, which can be thought of as an index of
refraction. The function $g : \Gamma \to \mathbb{R}$ specifies the
boundary conditions, and is subject to certain compatibility
conditions~\cite{Bornemann:2006aa}.

One way of arriving at the eikonal equation is by approximating the
solution $u$ of the \emph{Helmholtz equation}:
\begin{equation}
  \Big(\Delta + \omega^2 s(\m{x})^2\Big) u(\m{x}) = 0,
\end{equation}
with the WKB ansatz:
\begin{equation}
  u(\m{x}) \sim \alpha(\m{x}) e^{i \omega \tau(\m{x})},
\end{equation}
where $\omega$ is the frequency~\cite{Popov:2002aa}. As
$\omega \to \infty$, this asymptotic approximation is $O(\omega^{-1})$
accurate. This is referred to as the \emph{geometric optics}
approximation~\cite{Benamou:2003aa}. The level sets of $\tau$ denote
the arrival times of bundles of rays, and the amplitude $\alpha$, which
satisfies the transport equation:
\begin{equation}\label{eq:transport-equation}
  \alpha(\m{x}) \Delta \tau(\m{x}) + 2 \nabla \tau(\m{x})^\top \nabla \alpha(\m{x}) = 0,
\end{equation}
describes the attenuation of the amplitude of the wavefront due to the
propagation and geometric spreading of rays. The characteristics
(\emph{rays}) of the eikonal equation satisfy the raytracing ODEs.

The solution of the eikonal equation is given by Fermat's principle:
\begin{equation}\label{eq:fermats-principle}
  \tau(\m{x}) = \min_{\substack{\m{y} \in \Gamma \\ \mpsi : [0, 1] \to \Omega \\ \mpsi(0) = \m{y}, \mpsi(1) = \m{x}}} \left\{\tau(\m{y}) + \int_0^1 s(\mpsi(\sigma)) \|\mpsi'(\sigma)\|d\sigma\right\}.
\end{equation}
Observe that this equation is recursive, suggesting a connection with
dynamic programming and Bellman's principle of optimality. Indeed, the
path $\mpsi$ is a ray whose Lagrangian and Hamiltonian are:
\begin{equation}
  \mathcal{H}(\m{x}, \nabla \tau(\m{x})) = \frac{\|\nabla\tau(\m{x})\|^2 - s(\m{x})^2}{2} = 0, \qquad \mathcal{L}(\m{x}, \dot{\m{x}}) = s(\m{x})\|\dot{\m{x}}\|.
\end{equation}
This provides the connection between the Eulerian perspective given by
the eikonal equation, Fermat's principle, and the Lagrangian view
provided by raytracing.

\section{Related work}

The quintessential numerical method for solving the eikonal equation
is the \emph{fast marching method}~\cite{Sethian:1996aa}. We
discretize $\Omega$ into a grid of nodes $\Omega_h$, where $h > 0$ is
the characteristic length scale of elements in $\Omega_h$. Let
$T : \Omega_h \to \mathbb{R}$ be the numerical eikonal. To compute
$T$, equation \eqref{eq:eikonal-equation} is discretized using
first-order finite differences and the order in which individual
values of $T$ are relaxed is determined using a variation of
Dijkstra's algorithm for solving the single source shortest paths
problem~\cite{Sethian:1999ac,Sethian:1999ab}. If $N = |\Omega_h|$,
then the fast marching method solves \eqref{eq:eikonal-equation} in
$O(N \log N)$ with $O(h \log \tfrac{1}{h})$ worst-case
accuracy~\cite{Zhao:2005aa}. The logarithmic factor only appears when
rarefaction fans are present: e.g., point source boundary data, or if
the wavefront diffracts around a singular corner or edge. In these
cases, full $O(h)$ accuracy can be recovered by proper initialization
near rarefaction fans, or by employing a variety of factoring
schemes~\cite{Fomel:2009aa,Luo:2014aa,Qi:2019aa}.

It is also possible to solve the eikonal equation using
semi-Lagrangian numerical methods, in which the ansatz
\eqref{eq:fermats-principle} is discretized and applied
locally~\cite{Tsitsiklis:1995aa,Potter:2019ab}. For instance, at a
point $\hat{\m{x}} \in \Omega_h$, we consider a neighborhood of points
$\nb(\m{x}) \subseteq \Omega_h$, assume that $\tau$ is fixed over the
``surface'' of this neighborhood, and approximate
\eqref{eq:fermats-principle}. As an example, if $\nb(\hat{\m{x}})$
consists of its $2n$ nearest neighbors, if we linearly interpolate
$\tau$ over the facets of $\conv(\nb(\xhat))$, and discretize the
integral in \eqref{eq:fermats-principle} using a right-hand rule, the
resulting solver is equivalent to the fast marching
method~\cite{Sethian:2003aa}.

The Eulerian approach has generally been favored when developing
higher-order solvers for the eikonal equation~\cite{Zhang:2006aa}. The
eikonal equation is discretized using higher-order finite difference
schemes and solved in the same manner as the fast marching method or
using a variety of appropriate iterative schemes. Unfortunately, these
approaches presuppose a regular grid and require wide stencils.

Our goal is to develop solvers for the eikonal equation that are
\emph{high order}, are \emph{optimally local} (only use information
from the nodes in $\nb(\hat{\m{x}})$ to update $\hat{\m{x}}$), and are
flexible enough to work on \emph{unstructured meshes}. Using a
semi-Lagrangian approach based on a high-order discretization of
\eqref{eq:fermats-principle} allows us to do this.

This work was inspired by several lines of research. First, are
\emph{gradient-augmented level set methods} (or \emph{jet
  schemes})~\cite{Nave:2010aa,Seibold:2011aa}. Although developed for
solving time-dependent advection problems, trying to map ideas from
the \emph{time-dependent} to \emph{static} setting is natural, and
presented an intriguing challenge. Second, the idea of using a
semi-Lagrangian solver to construct a finite element solution to the
eikonal equation incrementally was
informative~\cite{Bornemann:2006aa}; while the authors only
constructed a first-order finte element approximation, attempting to
push past this formulation to obtain a higher-order solver is a
natural extension. Third, we were motivated by Chopp's idea of
building up piecewise bicubic interpolants locally while marching the
eikonal~\cite{Chopp:2001aa}; indeed, Chopp's work is mentioned in the
original work on jet schemes in a similar capacity.

\section{The jet marching method}

Label-setting algorithms~\cite{Chacon:2012aa}, such as the fast
marching method, compute an approximation to $\tau$ by marching a
numerical approximation $T : \Omega_h \to \mathbb{R}^n$ throughout the
domain. The boundary data $g$ is not always specified at the nodes of
$\Omega_h$. Let $\Gamma_h$ be a discrete approximation of
$\Gamma$. Once $T$ is computed at $\Gamma_h \subseteq \Omega_h$ with
sufficiently high accuracy, the solver begins to operate. To drive the
solver, a set of states $\{\far,\trial,\valid\}$ is used for
bookkeeping. We initially set:
\begin{equation}
  \state(\m{x}) = \begin{cases}
    \trial, & \mbox{if } \m{x} \in \Gamma_h, \\
    \far, & \mbox{otherwise.}
  \end{cases}
\end{equation}
The $\trial$ nodes are typically sorted by their $T$ value into an
array-based binary heap implementing a priority queue, although
alternatives have been explored~\cite{Gomez:2019aa}. At each step of
the iteration, the node $\m{x}$ with the minimum $T$ value is removed
from the heap, $\state(\m{x})$ is set to $\valid$, the $\far$ nodes in
$\nb(\m{x})$ have their state set to $\trial$, and each $\trial$ node
in $\nb(\m{x})$ is subsequently updated. We have additionally provided
a video online which shows the algorithm
running~\cite{Cameron:2020aa}.

From this, we can see that the value $T(\m{x})$ depends on the values
of $T$ at the nodes of a directed graph leading from $\m{x}$ back to
$\Gamma_h$, noting that $T(\m{x})$ can---and in general does---depend
on multiple nodes in $\nb(\m{x})$. This means that the error in $T$
accumulates as the solution propagates downwind from $\Gamma_h$. We
generally assume that the depth of the directed graph of updates
connecting each $\m{x} \in \Omega_h$ to $\Gamma_h$ is $O(h^{-1})$. The
error due to each update comes from two sources: the running error
accumulated in $T$, and the error incurred by approximating the
integral in \eqref{eq:fermats-principle}. For this reason, we would
expect the order of the global error of the solver to be one less than
the local error. However, the situation is more delicate because of
the complicated manner in which the errors mix (see
\Cref{fig:jmm4-pointwise-convergence}). Our numerical results in
\Cref{sec:numerical-results} indicate that $T$ converges with between
$O(h^2)$ and $O(h^3)$ accuracy, and $\nabla T$ converges with anywhere
between $O(h)$ and $O(h^3)$ accuracy.

Regardless, we assume that we only know the values of the eikonal and
some of its derivatives at the nodes $\m{x} \in \Omega_h$. To obtain
higher-order accuracy locally, we make use of piecewise Hermite
elements. In particular, at each node $\m{x}$, we approximate the
\emph{jet} of the eikonal; i.e., $\tau$ and a number of its
derivatives~\cite{Shilov:1971aa}. If we compute the jet with
sufficiently high accuracy when we set $\state(\m{x}) \gets \valid$,
we will be in a position to approximate $\tau$ using Hermite
interpolation locally over $\conv(\m{x}_1, \hdots, \m{x}_n)$. We
consider several variations on this idea.

\subsection{The general cost function}

Fix a point $\hat{\m{x}} \in \Omega_h$, thinking of it as the
\emph{update point}. To compute $T(\hat{\m{x}})$, we consider sets of
$\valid$ nodes
$\{\m{x}_1, \hdots, \m{x}_d\} \subseteq \nb(\hat{\m{x}})$, where
$1 \leq d \leq n$. The tuple of nodes
$(\hat{\m{x}}, \m{x}_1, \hdots, \m{x}_d)$ is an \emph{update} of
dimension $d$, and the collection of updates a \emph{stencil}. We
refer to the nodes $\{\m{x}_1, \hdots, \m{x}_d\}$ as the vertices of
the \emph{base of the update}. In some cases, such as on an
unstructured mesh, stencils may vary with $\xhat$.  Sufficient
conditions for the updates and stencils to be \emph{monotone causal}
have been studied~\cite{Kimmel:1998aa}. In particular, the cone
spanned by $\{\m{x}_1 - \xhat, \hdots, \m{x}_n - \xhat\}$ should fit
inside the nonnegative orthant after being
rotated~\cite{Sethian:2000aa,Sethian:2003aa}. This is easily satisfied
on a regular grid. For $O(h)$ solvers that do not make use of gradient
information, a variety of choices of stencils are monotone causal.

To describe a general update, without loss of generality we assume
$d = n$, and assume that the update nodes are in general
position. That is, if we choose $n$ nodes from
$\{\hat{\m{x}}, \m{x}_1, \hdots, \m{x}_n\}$, the remaining node does
not lie in their convex hull. We assume that we have access to a
sufficiently accurate approximation of $\tau$ over
$\conv(\m{x}_1, \hdots, \m{x}_n)$, call it $\mathsf{T}$. We
distinguish between $\mathsf{T}$ and $T$ in the following way:
$\mathsf{T}$ denotes the local numerical approximation of $\tau$ used
by a particular update, while $T$ denotes the global numerical
approximation of $\tau$. The two may not be equal to each
other. Indeed, $T$ is in general only defined on $\Omega_h$, while
$\mathsf{T}$ is only defined on $\conv(\m{x}_1, \hdots, \m{x}_n)$.

Let $\m{x}_{\m{\lambda}} \in \conv(\m{x}_1, \hdots, \m{x}_n)$, and let
$L = L_{\m{\lambda}} = \|\hat{\m{x}} - \m{x}_{\m{\lambda}}\|$. Recall
that $\mpsi : [0, L] \to \Omega$ is the curve
minimizing~\eqref{eq:fermats-principle} for a particular choice of
$\xlam$. We approximate $\mpsi$ with a cubic parametric curve
$\mphi : [0, L] \to \Omega$ such that:
\begin{equation}
  \mphi(0) = \xlam, \qquad \mphi(L) = \xhat, \qquad \mphi'(0) \sim \tlam, \qquad \mphi'(L) \sim \that,
\end{equation}
and where $\tlam$ and $\that$ are tangent vectors which enter as
parameters. Note that $\mphi'(0)$ and $\mphi'(L)$ may not be exactly
equal to $\tlam$ and $\that$.

We approximate the integral in \eqref{eq:fermats-principle} over
$\mphi$ using Simpson's rule. This gives the cost functional:
\begin{equation}\label{eq:general-cost-function}
  F(\mphi) = \mathsf{T}(\xlam) + \frac{L}{6}\Big[s(\xlam)\|\mphi'(0)\| + 4 s\big(\mphi_{\half}\big)\big\|\mphi'_{\half}\big\| + s(\xhat)\|\mphi'(L)\|\Big],
\end{equation}
where $\mphi{}_{\half} = \mphi(L/2)$ and
$\mphi'\hspace{-0.2em}{}_{\half} = \mphi'(L/2)$.  We have not yet made this
well-defined. To do so, we must specify $\mathsf{T}$, $\tlam$, and
$\that$. We describe several different ways of doing this in the
following sections.

\subsection{Computing $\nabla T(\xhat)$}

A minimizing extremal $\mpsi$ of Fermat's integral is a characteristic
of the eikonal equation. A simple but important consequence of this is
that its tangent vector is locally parallel to $\nabla \tau$. Hence:
\begin{equation}
  s(\mpsi(\sigma)) \frac{\mpsi'(\sigma)}{\|\mpsi'(\sigma)\|} = \nabla \tau(\mpsi(\sigma)).
\end{equation}
After minimizing $F$, we will have found an optimal value of
$\that$. We can then set:
\begin{equation}
  \nabla T(\xhat) \gets s(\xhat) \that.
\end{equation}
This lets us march the gradient of the eikonal locally along with the
eikonal itself.

\subsection{Parametrizing $\mphi$}

\begin{figure}
  \centering
  \includegraphics{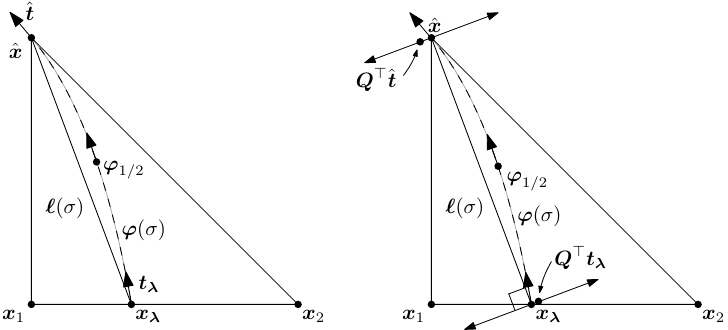}
  \caption{Two approaches to parametrizing a cubic curve approximating
    the characteristic $\mphi$ leading from $\xlam$ to $\xhat$ when
    numerically minimizing Fermat's integral to compute $T(\xhat)$ and
    $\nabla T(\xhat)$. \emph{Left}: $\mphi$ is a cubic parametric
    curve with boundary data set directly from $\xlam, \xhat, \tlam$,
    and $\that$. \emph{Right}: $\mphi$ is the graph of a function in
    the orthogonal complement of
    $\operatorname{range}(\mell')$.}\label{fig:parametrizations}
\end{figure}

We consider two methods of choosing $\mphi$ (see
Figure~\ref{fig:parametrizations}). First, let $\mell$ be interval
connecting $\xlam$ to $\xhat$, parametrized by arc length, and define:
\begin{equation}
  \mell(\sigma) = \xlam + \sigma \mell', \qquad \mell' = (\xhat - \xlam)/L_{\lam}.
\end{equation}

\paragraph{Using a cubic parametric curve}
For one approach, we define:
\begin{equation}
  \mphi(\sigma) = \mell(\sigma) + \dphi(\sigma),
\end{equation}
where $\dphi : [0, L] \to \Omega$ is a perturbation away from
$\mell$ that satisfies:
\begin{equation}
  \dphi(0) = 0, \quad \dphi(L) = 0, \quad \dphi'(0) = \tlam - \mell', \quad \dphi'(L) = \that - \mell'.
\end{equation}
We can explicitly write $\dphi$ as:
\begin{equation}
  \dphi(\sigma) = \big(\tlam - \mell'\big) K_0(\sigma) + \big(\that - \mell'\big) K_1(\sigma),
\end{equation}
where $K_0, K_1 : [0, L] \to \mathbb{R}$ are Hermite basis functions
such that:
\begin{equation}
  \begin{split}
    K_0(0) = 0 = K_0(L), \qquad &K_1(0) = 0 = K_1(L), \\
    K_0'(0) = 1 = K_1'(L), \qquad &K_1'(0) = 0 = K_0'(L).
  \end{split}
\end{equation}
Explicitly, these are given by:
\begin{equation}\label{eq:K-formulas}
  K_0(\sigma) = \sigma - 2 \frac{\sigma^2}{L} + \frac{\sigma^3}{L^2}, \qquad K_1(\sigma) = \frac{-\sigma^2}{L} + \frac{\sigma^3}{L^2}.
\end{equation}

Let $\tlam, \that \in \mathbb{S}^{n-1}$ so that
$\|\tlam\| = 1 = \|\that\|$.  As $L \to 0$, this results in a curve
that is approximately parametrized by arc length: i.e.,
$\|\mphi'(\sigma)\| \to 1$ for all $\sigma$ such that
$0 \leq \sigma \leq L$~\cite{Floater:2006aa}. This simplifies the
general cost function given by~\eqref{eq:general-cost-function} to:
\begin{equation}\label{eq:F-curve-form}
  F(\mphi) = \mathsf{T}(\xlam) + \frac{L}{6} \Big[s(\xlam) + 4 s\big(\mphi_{\half}\big)\big\|\mphi'_{\half}\big\| + s(\xhat)\Big].
\end{equation}
Using~\eqref{eq:K-formulas}, $\mphi_{\half}$ and $\mphi'_{\half}$ can
be written:
\begin{equation}
  \mphi_{\half} = \frac{\xlam + \xhat}{2} + \frac{L}{8}\big(\tlam - \that\big), \qquad \mphi'_{\half} = \frac{3}{2} \mell' - \frac{\tlam + \that}{4}.
\end{equation}
Note that $\mphi{}_{\half} \sim (\xlam + \xhat)/2$ and
$\mphi'\hspace{-0.2em}{}_{\half} \sim \mell'$ as $L \to 0$ if we
assume that the wavefront is well-approximated by a plane wave near
the update, since in this case $\tlam \sim \that \sim \mell'$.

\paragraph{Parametrizing $\mphi$ as the graph of a function} We can
also define the perturbation away from $\mell$ as the graph of a
function; i.e., we assume that the perturbation is orthogonal to
$\mell'$. Letting $\m{Q} \in \mathbb{R}^{n \times (n-1)}$ be an
orthogonal matrix such that $\m{Q}^\top \mell' = 0$, and letting
$\m{\zeta} : [0, L] \to \mathbb{R}^{n-1}$ be a curve specifying the
components of the perturbation in this basis, we choose
$\dphi(\sigma) = \m{Q}\m{\zeta}(\sigma)$ so that:
\begin{equation}
  \mphi(\sigma) = \mell(\sigma) + \m{Q} \m{\zeta}(\sigma).
\end{equation}
where $\m{\zeta}(\sigma) = \m{b}_0 K_0(\sigma) + \m{b}_1 K_1(\sigma)$.
In this approach, instead of $\that$ and $\tlam$, we optimize over
$\m{b}_0, \m{b}_1 \in \mathbb{R}^{n - 1}$. Now, noting that:
\begin{equation}
  \|\mphi'(\sigma)\| = \sqrt{\|\mell'\|^2 + \|\m{Q}\m{\zeta}'(\sigma)\|^2} = \sqrt{1 + \|\m{\zeta}'(\sigma)\|^2},
\end{equation}
we can write the cost functional $F$ as:
\begin{equation}\label{eq:F-graph-form}
  \begin{split}
    F(\mphi) &= \mathsf{T}(\xlam) + \frac{L}{6}\Bigg[s(\xlam)\sqrt{1 + \|\m{b}_0\|^2} \\
    &\hspace{3em} + 4 s(\mphi_{\half}) \sqrt{1 + \|{(\m{b}_0 +
        \m{b}_1)}/4\|^2} + s(\xhat)\sqrt{1 + \|\m{b}_1\|^2}\Bigg].
  \end{split}
\end{equation}

\paragraph{Trade-offs between the two parametrizations of $\mphi$}
When $\mphi$ is a cubic parametric curve, we run into an interesting
problem described in more detail by Floater~\cite{Floater:2006aa}. In
particular, the order of accuracy of $\mphi$ in approximating $\mpsi$
is limited by our parametrization of $\mphi$. If we parametrize
$\mphi$ over $\sigma \in [0, 1]$ (that is \emph{uniformly}), then the
interpolant is at most $O(h^2)$ accurate. If we parametrize it using a
\emph{chordal parametrization}, i.e. $\sigma \in [0, L]$, then it is
at most $O(h^4)$ accurate. Indeed, any Hermite spline using a chordal
parametrization over each of its segment is at most $O(h^4)$ accurate
globally. To design a higher order solver than this requires us to
parametrize $\mphi$ using a more accurate approximation of the arc
length of $\mpsi$ (consider, e.g., using a quintic parametric
curve). On the other hand, if we parametrize $\mphi$ as the graph of a
function, we can directly apply Hermite interpolation
theory~\cite{Stoer:2013aa}, and there is no such obstacle.

\section{Different types of minimization
  problems}\label{sec:minimization-problems}

In this section, we consider four different ways of using $F$ to pose
a minimization problem which would allow us to compute $T(\xhat)$. We
note that each of these formulations is compatible with the version of
$F$ where we take $\mphi$ to be a parametric curve \emph{and} where we
define it as the graph of a function orthogonal to
$\mell(\sigma)$. Altogether, this leads to eight different JMMs.

\subsection{Determining $\tlam$ by minimizing Fermat's
  integral}\label{ssec:full-minimization}

Since the optimal $\mphi$ is a characteristic of the eikonal equation,
one approach to setting $\tlam$ and $\that$ is to simply let them
enter into the cost function as free parameters to be optimized
over. This leads to the optimization problem:
\begin{equation}\label{eq:minimization-problem-1}
  \begin{split}
    \mbox{minimize} &\quad F(\xlam, \tlam, \that) \\
    \mbox{subject to} &\quad \xlam \in \conv(\m{x}_1, \hdots, \m{x}_n), \\
    &\quad \tlam, \that \in \mathbb{S}^{n-1},
  \end{split}
\end{equation}
if we parametrize $\mphi$ as a curve; or, if we parametrize $\mphi$ as
the graph of a function:
\begin{equation}\label{eq:minimization-problem-1-b}
  \begin{split}
    \mbox{minimize} &\quad F(\xlam, \m{b}_0, \m{b}_1) \\
    \mbox{subject to} &\quad \xlam \in \conv(\m{x}_1, \hdots, \m{x}_n), \\
    &\quad \m{b}_0, \m{b}_1 \in \mathbb{R}^{n-1},
  \end{split}
\end{equation}
For a $d$-dimensional update, the domain of each of these minimization
problems has dimension $(d - 1)(n - 1)^2$, since
$\dim(\conv(\m{x}_1, \hdots, \m{x}_d)) = d - 1$.

\subsection{Determining $\tlam$ from the eikonal
  equation}\label{ssec:eikonal-method}

When we compute updates, we only require high-order accurate jets over
$\conv(\m{x}_1, \hdots, \m{x}_n)$. This is a subset of $\Omega$ of
codimension one: an interval in 2D, or triangle in 3D. If we know $T$
and $\nabla T$ at the vertices of this set, then we can use Hermite
interpolation to compute $\mathsf{T}$. Unfortunately, this means that
we can only approximate directional derivatives of $T$ in the linear
span of this set. To compute $\nabla\mathsf{T}$, we need to recover
the directional derivative normal to the facet.

Let $\tilde{\m{V}} \in \mathbb{R}^{n \times (n-1)}$ be an orthogonal matrix such
that:
\begin{equation}
  \range(\tilde{\m{V}}) = \range\hspace{-0.1em}\left(\begin{bmatrix} \m{x}_2 - \m{x}_1 & \cdots & \m{x}_n - \m{x}_1 \end{bmatrix}\right),
\end{equation}
and let $\m{v} \in \mathbb{R}^n$ be a unit vector such that
$\tilde{\m{V}}{}^\top \m{v} = 0$. Let $\nabla_{\tilde{\m{V}}}$ be the
gradient restricted to the range of $\tilde{\m{V}}$, and likewise let
$d_{\m{v}}$ denote the $\m{v}$-directional derivative. Then, from the
eikonal equation, we have:
\begin{equation}
  s(\m{x})^2 = \|\nabla\tau(\m{x})\|^2 = |d_{\m{v}}\tau(\m{x})|^2 + \|\nabla_{\tilde{\m{V}}}\tau(\m{x})\|^2.
\end{equation}
To recover $\nabla \tau(\m{x})$, first note that $\nabla \tau(\m{x})$
should point in the same direction as $\mell'$. Choosing $\m{v}$ so
that $\m{v}^\top \mell' > 0$, we get:
\begin{equation}\label{eq:v-directional-derivative}
  d_{\m{v}}\tau(\m{x}) = \sqrt{s(\m{x})^2 - \|\nabla_{\tilde{\m{V}}} \tau(\m{x})\|^2}.
\end{equation}
Letting $\m{V} = \begin{bmatrix} \m{v} & \tilde{\m{V}} \end{bmatrix}$,
equation \eqref{eq:v-directional-derivative} combined with
$\nabla \tau(\m{x}) = \m{V} \nabla_{\m{V}} \tau(\m{x})$ gives us a
means of recovering $\nabla\tau(\m{x})$ from
$\nabla_{\tilde{\m{V}}} \tau(\m{x})$.

Using this technique, we can pose the following optimization problem:
\begin{equation}\label{eq:minimization-problem-2}
  \begin{split}
    \mbox{minimize} &\quad F(\xlam, \that) \\
    \mbox{subject to} &\quad \xlam \in \conv(\m{x}_1, \hdots, \m{x}_n), \\
    &\quad \that \in \mathbb{S}^{n-1},
  \end{split}
\end{equation}
or, optimizing over $\m{b}_1$ directly:
\begin{equation}\label{eq:minimization-problem-2-b}
  \begin{split}
    \mbox{minimize} &\quad F(\xlam, \m{b}_1) \\
    \mbox{subject to} &\quad \xlam \in \conv(\m{x}_1, \hdots, \m{x}_n), \\
    &\quad \m{b}_1 \in \mathbb{R}^{n-1}
  \end{split}
\end{equation}
For each $\xlam$, we set:
\begin{equation}
  \tlam \gets \frac{\m{V} \nabla_{\m{V}} \mathsf{T}(\xlam)}{\left\|\m{V} \nabla_{\m{V}} \mathsf{T}(\xlam)\right\|}.
\end{equation}
Note that the dimension of a $d$-dimensional update based on this
minimization problem is $(d - 1)(n - 1)$

\subsection{Determining $\tlam$ by marching cell-based
  interpolants}\label{ssec:cell-based}

Another approach is to march cells that approximate the jet of the
eikonal at each point. For example, if we have constructed a finite
element interpolant using $\valid$ data on a cell whose boundary
contains $\conv(\m{x}_1, \hdots, \m{x}_n)$ then we can evaluate its
gradient to obtain:
\begin{equation}
  \tlam \gets \frac{\nabla \mathsf{T}(\xlam)}{\left\|\nabla\mathsf{T}(\xlam)\right\|}.
\end{equation}
We can combine this approach with the cost functional given by
\eqref{eq:minimization-problem-2}, albeit with a modified $\tlam$. We
elaborate on how we march cells in section~\ref{sec:cell-marching}. An
advantage of this approach is that it allows one to simultaneously
march the second partials of $T$.

\subsection{A simplified method using a quadratic curve}\label{ssec:quadratic}

In some cases, in particular if the speed of sound is linear, i.e.:
\begin{equation}
  c(\m{x}) = \frac{1}{s(\m{x})} = c_0 + \m{c}^\top \m{x}, \qquad c_0 \in \mathbb{R}, \qquad \m{c} \in \mathbb{R}^n,
\end{equation}
the characteristic $\mpsi$ is well-approximated by a quadratic. In
this case, we again have a cost functional of the form
\eqref{eq:minimization-problem-2}.

If $\mphi$ is parametrized as curve, we set $\tlam$ to be the
reflection of $\that$ across $\mell'$:
\begin{equation}
  \tlam = -\big(\m{I} - 2 \hspace{0.1em} \mell' \mell'^\top\big) \that,
\end{equation}
giving $\tlam + \that = 2 \mell'\mell'^\top \that$ and
$\tlam - \that = -2\big(\m{I} - \mell'\mell'^\top\big)\that$. Then:
\begin{equation}
  \mphi_{\half} = \frac{\xlam + \xhat}{2} - \frac{L}{4}\big(\m{I} - \mell' \mell'^\top\big) \that, \qquad \mphi'_{\half} = \frac{3 + \mell'^\top \that}{2} \mell'.
\end{equation}
This simplifies $F$ given by~\eqref{eq:F-curve-form} to:
\begin{equation}
  F(\xlam, \that) = \mathsf{T}(\xlam) + \frac{L}{6}\bigg[s(\xlam) + 2\big(3 + \mell'^\top \that\big) s\Big(\frac{\xlam + \xhat}{2} - \frac{L}{4}\big(\m{I} - \mell' \mell'^\top\big)\Big) + s(\xhat)\bigg].
\end{equation}

If $\mphi$ is parametrized as the graph of a function, then:
\begin{equation}
  \m{\zeta}_{\half} = \frac{\xhat + \xlam}{2} + \frac{L}{4} \m{Q}^\top \that, \qquad \m{\zeta}' = 0,
\end{equation}
simplifying the version of $F$ in~\eqref{eq:F-graph-form} to:
\begin{equation}
  F(\xlam, \that) = \mathsf{T}(\xlam) + \frac{L}{6}
  \bigg[\big(s(\xlam) + s(\xhat)\big) \sqrt{1 +
    \|\m{b}_0\|^2} + 4 s(\mphi_{\half})\bigg].
\end{equation}
since $\m{Q}^\top\that = \m{b}_0 = -\m{b}_1$.

\subsection{Other approaches}

We tried two other approaches which failed to provide satisfactory
results:
\begin{itemize}
\item A combination of the quadratic simplification in
  subsection~\ref{ssec:quadratic} with the methods in
  subsections~\ref{ssec:eikonal-method} or~\ref{ssec:cell-based}. In
  this case, we use our knowledge of $\nabla T(\xlam)$ along the base
  of the update to choose $\that$ and $\tlam$. This reduces the
  dimensionality of the cost function to $d - 1$. However, except for
  in special cases (e.g.\ $s \equiv 1$), this propagates errors in a
  manner that causes the solver to diverge; or, at best, allows it
  converge with $O(h)$ accuracy. We note that if $s \equiv 1$, still
  simpler methods can be used, so this combination of approaches does
  not seem to be useful.
\item We can extract not only $\tlam$ from the Hermite interpolant on
  $\conv(\m{x}_1, \hdots, \m{x}_n)$, but also $\mphi''(0)$. From the
  Euler-Lagrange equations for the eikonal equation, we obtain:
  \begin{equation}
    \m{Q}^\top \nabla s(\xlam) = s(\xlam) \frac{\m{Q}^\top \mphi''(0)}{\|\mphi'(0)\|}.
  \end{equation}
  With $\mphi$ parametrized as a graph, we have
  $\m{Q}^\top \mphi''(0) = \m{\zeta}''(0)$, giving:
  \begin{equation}
    \m{\zeta}''(0) = \frac{\m{Q}^\top \nabla s(\xlam) {(1 + \|\m{Q}^\top \xlam\|^2)}}{s(\xlam)}.
  \end{equation}
  This completely defines $\mphi$ as the graph of a cubic polynomial
  using the graph parametrization. Unfortunately, this method diverges
  for the same reason as the method described in the previous bullet.
\end{itemize}

\subsection{A warm start}\label{sec:warm-start}

Certain of the optimization problems in the preceding section are
clearly nonconvex; e.g., \eqref{eq:minimization-problem-1} is
nonconvex since its domain, the product of
$\conv(\m{x}_1, \hdots, \m{x}_n)$ and two copies of
$\mathbb{S}^{n-1}$, is nonconvex. For $h$ small, if $\xlam$ is nearly
optimal, then the optimal local ray $\mphi$ should be close to
$\mell$, the straight line segment connecting $\xlam$ and
$\xhat$. This suggests an approach to finding an initial iterate (a
warm start) for \eqref{eq:minimization-problem-1} or the other
optimization problems considered (i.e.,
\eqref{eq:minimization-problem-1-b},
\eqref{eq:minimization-problem-2}, and
\eqref{eq:minimization-problem-2-b}). In 2D, following a similar
approach to the simplified midpoint rule (denoted ``\texttt{mp0}'')
rule used in our earlier work on ordered line integral methods (OLIMs)
for the solving the eikonal equation~\cite{Potter:2019ab}, we let
$\mathsf{T}$ be a cubic Hermite polynomial approximation of $\tau$
over $\lambda \in [0, 1]$ and approximate \eqref{eq:fermats-principle}
by solving:
\begin{equation}\label{eq:warm-start-minimization-problem}
  \begin{split}
    \mbox{minimize} &\quad \mathsf{T}(\lambda) + \frac{L}{2} \Big(s(\m{x}_\lambda) + s(\xhat)\Big), \\
    \mbox{subject to} &\quad 0 \leq \lambda \leq 1.
  \end{split}
\end{equation}
After solving this problem, we can compute an initial iterate for
\eqref{eq:minimization-problem-1} from $\lambda^*$, the optimum of
\eqref{eq:warm-start-minimization-problem}. For example, we set
$\m{x}_{\lambda} \gets \m{x}_{\lambda^*}$; we set $\m{t}_{\lambda}$ using $\lambda^*$ and
one of the approaches outlined in the preceding sections; and, if
required, we set
$\that \gets (\xhat - \m{x}_{\lambda^*})/\|\xhat -
\m{x}_{\lambda^*}\|$. In practice,
\eqref{eq:warm-start-minimization-problem} can be solved rapidly and
robustly using a rootfinder.

\subsection{Optimization algorithms}

We do not dwell on the details of how to numerically solve the
minimization problems in the preceding sections. We make some general
observations:
\begin{itemize}
\item These optimization problems are very easy to solve---what's
  costly is that we have to solve $O(N)$ of them. As $h \to 0$, they
  are strictly convex and well-behaved. Empirically, Newton's method
  converges in $O(1)$ steps (typically fewer than 5 with a well-chosen
  warm start---see section~\ref{sec:warm-start}). We leave a detailed
  comparison of different approaches to numerically solving these
  optimization problems for future work.
\item The gradients and Hessians of these cost functions are somewhat
  complicated. Programming them can be tricky and tedious, suggesting
  that automatic differentiation may be a worthwhile
  approach~\cite{Neidinger:2010aa,Griewank:2008aa}.
\item The constraint $\xlam \in \conv(\m{x}_1, \hdots, \m{x}_n)$
  corresponds to a set of linear inequality constraints, which are
  simple to incorporate. Because of the form of these constraints,
  checking the KKT conditions at the boundary is cheap and easy~\cite{Potter:2019ab,Yang:2019aa}. See
  the next section on skipping updates.
\item The constraints $\tlam, \that \in \mathbb{S}^{n-1}$ are
  nonlinear; however, they can be eliminated. If $n = 2$, then we can
  set $\that = (\cos(\theta), \sin(\theta))$, letting
  $\theta \in \mathbb{R}$. For $n > 2$, we can use a Riemannian
  Newton's method for minimization on $\mathbb{S}^{n-1}$, which is
  simple to implement and known to converge
  superlinearly~\cite{Absil:2009aa}. Alternatively, we could use
  spherical coordinates, although the expressions become unwieldy.
\end{itemize}

\section{Hierarchical update algorithms}

Away from shocks, where multiple wavefronts collide, exactly one
characteristic will pass through a point $\xhat$. When we minimize $F$
over each update in the stencil, the characteristic will pass through
the base of the minimizing update, or possibly through the boundary of
several adjacent updates. We can use this fact to sequence the updates
that are performed to design a work-efficient solver. In our previous
work on \emph{ordered line integral methods} (\emph{OLIMs}), we
explored variations of this idea~\cite{Potter:2019ab,Yang:2019aa}. An
approach that works well is the \emph{bottom-up} update algorithm.

\begin{figure}
  \centering
  \includegraphics{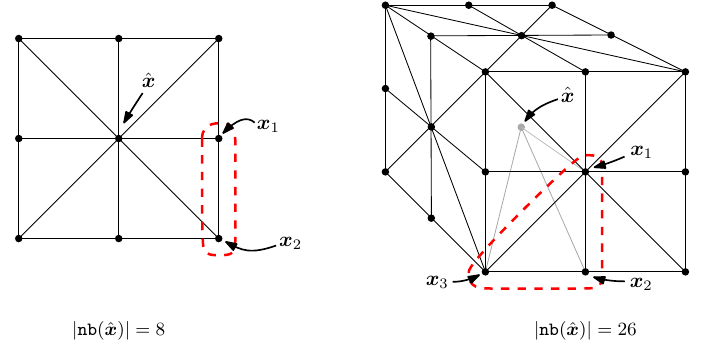}
  \caption{The neighborhoods typically used by semi-Lagrangian solvers
    in 2D and 3D on a regular grid. These are the stencils used by
    Tsitsiklis's algorithm and two of the OLIM
    stencils~\cite{Tsitsiklis:1995aa,Potter:2019ab}. \emph{Left}:
    ``\texttt{olim8}'' in $\mathbb{R}^2$. This is the 8-point stencil
    used in this paper. \emph{Right}: ``\texttt{olim26}'' in
    $\mathbb{R}^3$.}\label{fig:neighborhoods}
\end{figure}

To fix the idea in 3D, consider $\nb(\m{x})$ as shown in Figure\
\ref{fig:neighborhoods}, for which $|\nb(\m{x})| = 26$. There are 26
``line'' updates, where $d = 1$. To start with, each $\valid$ line
update is done, and $\m{x}_1$ for the minimizing line update is
recorded. Next, we fix $\m{x}_1$ and perform ``triangle'' updates
($d = 2$) where $\m{x}_2$ is varying. In this case, we can restrict
the number of triangle updates that are done by assuming either that
$(\m{x}_1, \m{x}_2)$ is an edge of mesh discretizing the surface of
the 3D stencil shown in Figure\ \ref{fig:neighborhoods}, or that
$\|\m{x}_1 - \m{x}_2\|$ is small enough (measuring the distance of
these two points in different norms leads to a different number of
triangle updates---we find the $\ell_1$ norm to work well). Finally,
we fix $\m{x}_2$ corresponding to the minimizing triangle update, and
do tetrahedron updates containing $\m{x}_1$ and $\m{x}_2$. Throughout
this process, $\m{x}_1, \m{x}_2$, and $\m{x}_3$ must all be $\valid$.

We emphasize that our work-efficient OLIM update algorithms work
equally well for the class of algorithms developed here. The main
differences between the JMMs studied here and the earlier OLIMs are
the cost functionals and the we way approximate $T$.

\section{Initialization methods}

A common problem with the convergence of numerical methods for solving
the eikonal equation concerns how to treat rarefaction fans. Our
numerical tests consist of point source problems, around which a
rarefaction forms. A standard approach is to introduce the factored
eikonal equation~\cite{Fomel:2009aa,Luo:2014aa,Qi:2019aa}. If a point
source is located at $\m{x}^\circ \in \Omega_h$ and if we set
$\Gamma_h = \{\m{x}^\circ\}$, then we let
$d(\m{x}) = \|\m{x} - \m{x}^\circ\|$ and use the ansatz:
\begin{equation}
  \tau(\m{x}) = z(\m{x}) + d(\m{x}), \qquad \m{x} \in \Omega.
\end{equation}
We insert this into the eikonal equation, modifying our numerical
methods as necessary, and solve for $z(\m{x})$ instead. This is not
complicated---see our previous work on OLIMs for solving the eikonal
equation to see how the cost functions should generally be
modified~\cite{Potter:2019ab,Yang:2019aa}.

Yet another approach would be to solve the characteristic equations to
high-order for each $\m{x}$ in such a ball. This would require solving
$O(N)$ boundary value problems, each discretized into $O(N^{1/3})$
intervals, resulting in an $O(N^{4/3})$ cost overall (albeit with a
very small constant). One issue with this approach is that it only
works well if the ball surrounding $\m{x}^\circ$ is contained in the
interior of $\Omega$. For our numerical experiments, we simply
initialize $T$ and $\nabla T$ to the correct, ground truth values in a
ball or box of constant size centered at $\m{x}^\circ$.

\section{Cell marching}\label{sec:cell-marching}

Of particular interest is solving the transport equation governing the
amplitude $\alpha$ while simultaneously solving the eikonal
equation. Equation \eqref{eq:transport-equation} can be solved using
upwind finite differences~\cite{Benamou:1996aa} or paraxial
raytracing~\cite{Popov:2002aa}. We prefer the latter approach since it
can be done locally, using the characteristic path $\mphi$ recovered
when computing $T(\xhat)$. Either approach requires accurate second
derivative information (we need $\Delta T$ for upwind finite
differences, or $\nabla^2 T$ for paraxial raytracing).

For the purposes of explanation and our numerical tests, we consider a
rectilinear grid with square cells in $\mathbb{R}^2$. On each cell,
our goal is to build a bicubic interpolant, approximating
$T(\m{x})$. This requires knowing $T, \nabla T$, and $T_{xy}$ at each
cell corner. If we know these values with $O(h^{4-p})$ accuracy, where
$p$ is the order of the derivative, then the bicubic is $O(h^{4-p})$
accurate over the cell. So far, we have described an algorithm that
marches $T$ and $\nabla T$, which together constitute the \emph{total
  1-jet}. We now show how $T_{xy}$ can also be marched, allowing us to
march the \emph{partial 1-jet}.\footnote{The total $k$-jet of a
  function $f$ is the set $\{\partial^{\m{\alpha}}f\}_{\m{\alpha}}$,
  where $\|\m{\alpha}\|_1 \leq k$; the partial $k$-jet is
  $\{\partial^{\m{\alpha}} f\}_{\m{\alpha}}$ where
  $\|\m{\alpha}\|_\infty \leq k$.}

\begin{figure}
  \centering
  \includegraphics[width=\textwidth]{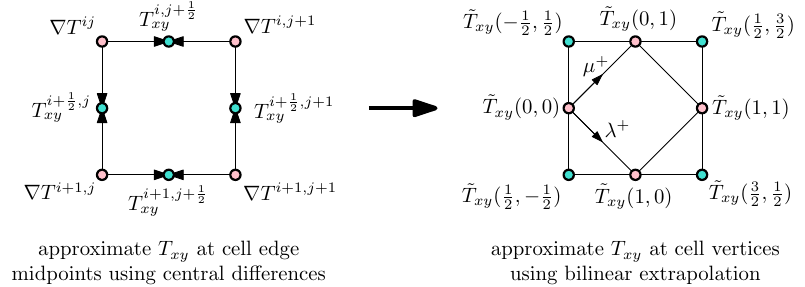}
  \vspace{-1em}
  \caption{\emph{Cell-based interpolation.} To approximate the mixed
    second partials of a function with $O(h^2)$ accuracy from $O(h^3)$
    accurate gradient values available at the corners of a cell, the
    following method of using central differences to approximate the
    mixed partials at the midpoints of the edges of the cell, followed
    by bilinear extrapolation, can be used.}\label{fig:estimate-Txy}
\end{figure}

Let $\m{x}_{ij}$ with $(i, j) \in \{0, 1\}^2$ denote the corners of a
square cell with sides of length $h$, and assume that we know
$\nabla T(\m{x}_{ij})$ with $O(h^3)$ accuracy. We can use the
following approach to estimate $T_{xy}(\m{x}_{ij})$ at each corner:
\begin{itemize}
\item First, at the midpoints of the edges oriented in the $x$
  direction (resp., $y$ direction), approximate $T_{xy}$ using the
  central differences involving $T_y$ (resp., $T_x$) at the
  endpoints. This approximation is $O(h^2)$ accurate at the midpoints.
\item Use bilinear extrapolation to reevaluate $T_{xy}$ at the corners
  of the cell, yielding $T_{xy}(\m{x}_{ij})$, also with $O(h^2)$
  accuracy.
\end{itemize}
This procedure is illustrated in Figure\ \ref{fig:estimate-Txy}.

\begin{figure}
  \centering
  \includegraphics[width=\linewidth]{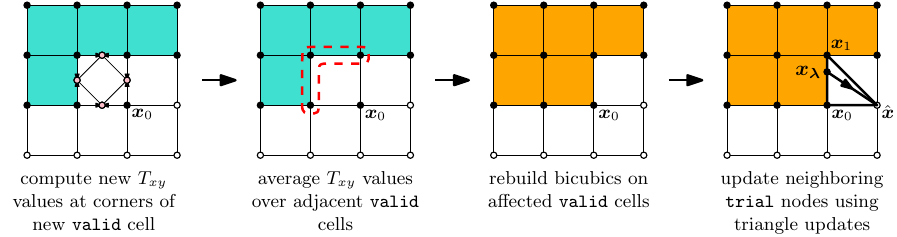}
  \vspace{-1em}
  \caption{\emph{Local cell marching.} After computing values of
    $T_{xy}$ as shown in Figure~\ref{fig:estimate-Txy} (\emph{left}),
    to ensure continuity of the global interpolant, nodal values
    incident on the newly \texttt{valid} cell (containing $x_0$) can
    be recomputed by averaging over $T_{xy}$ values taken from
    incident \texttt{valid} cells (\emph{middle}). Finally, a bicubic
    cell-based interpolant is constructed
    (\emph{right}).}\label{fig:rebuild-cells}
\end{figure}

One issue with this approach is that it results in a piecewise
interpolant that is only $C^1$ globally. That is, if we estimate the
value of $T_{xy}$ at a corner from each of the cells which are
incident upon it, we will get different values in general. To compute
a globally $C^2$ piecewise interpolant, we can average $T_{xy}$ values
over incident $\valid$ cells, where we define a $\valid$ cell to be a
cell whose vertices are all $\valid$. How to do this is shown in
Figure\ \ref{fig:rebuild-cells}.

The idea of approximating the partial 1-jet from the total 1-jet in an
optimally local fashion by combining central differences with bilinear
extrapolation, and averaging nodal values over adjacent cells to
increase the degree of continuity of the interpolant, is borrowed from
Seibold et al.~\cite{Seibold:2011aa}. However, applying this idea in
this context, and doing the averaging in an upwind fashion is novel.

The scheme arrived at in this way is no longer optimally
local. However, the sequence of operations described here can be done
on an unstructured triangle or tetrahedron mesh. This makes this
approach suitable for use with an unstructured mesh that conforms to a
complicated boundary. We should mention here that our approach to
estimating $T_{xy}$ is referred to as \emph{twist estimation} in the
\emph{computer-aided design} (CAD) community~\cite{Farin:2014aa},
where other approaches have been
proposed~\cite{Brunet:1985aa,Hagen:1987aa}. We leave adapting these
ideas to the present context for future work.

\subsection{Marching the amplitude} In this section we show how to
compute a numerical approximation of $\alpha$, denoted
$A : \Omega_h \to \mathbb{R}$. One simple approach would be to
discretize~\eqref{eq:transport-equation} using upwind finite
differences and compute $A(\xhat)$ using $\valid$ nodes after
$T(\xhat)$ and $\nabla T(\xhat)$. One potential shortcoming of this
approach is that $A$ is singular at caustics. Instead, we will explore
using paraxial raytracing to compute $A$ in this
section~\cite{Popov:1978aa}. The background material on paraxial
raytracing used in this section can be found in more detail in M.\
Popov's book~\cite{Popov:2002aa}.

The basic idea of paraxial raytracing is to consider a fixed, central
ray, which we denote $\mphi_0$, and a surrounding tube of rays,
parametrized by:
\begin{equation}
  \mphi(\sigma, \m{q}) = \mphi_0(\sigma) + \m{E}(\sigma)\m{q},
\end{equation}
where $\m{E} : [0, L] \to \mathbb{R}^{n \times (n-1)}$ is an
orthogonal matrix such that $\m{E}^\top \mphi_0' \equiv 0$. For each
$\m{q}$, the corresponding ray should satisfy the Euler-Lagrange
equations for~\eqref{eq:eikonal-equation}. If we let
$c_0(\sigma) = c(\mphi_0(\sigma))$, where $c = 1/s$, then $\m{q}$
along with the conjugate momenta $\m{p}$ (the exact form of which is
not important in this instance) will satisfy:
\begin{equation}\label{eq:hamiltons-equations-ray-tube}
  \begin{bmatrix}
    d\m{q}/d\sigma \\
    d\m{p}/d\sigma
  \end{bmatrix} = \begin{bmatrix}
    & c_0(\sigma) \m{I} \\
    \frac{-1}{c_0(\sigma)^2} \left. \frac{\partial^2 c}{\partial \m{q} \partial \m{q}^\top} \right|_{\m{q}=0} &
  \end{bmatrix} \begin{bmatrix}
    \m{q} \\ \m{p}
  \end{bmatrix}.
\end{equation}
If we let
$\m{Q}(\sigma), \m{P}(\sigma) : [0, L] \to \mathbb{R}^{(n-1) \times
  (n-1)}$ be a linearly independent set of solutions
to~\eqref{eq:hamiltons-equations-ray-tube}, then along the central
ray, the amplitude satisfies:
\begin{equation}\label{eq:amplitude-update}
  A(\mphi_0(\sigma)) = \sqrt{\frac{c_0(\sigma)}{|\det(\m{Q}(\sigma))|}} A(\mphi_0(0)).
\end{equation}
Note that when we compute an update, we obtain a cubic path $\mphi$
approximating a ray of~\eqref{eq:eikonal-equation}, such that
$\mphi(0) = \xlam$ and $\mphi(L) = \xhat$.

The quantity $|\det(\m{Q}(\sigma))|$ is known as the \emph{geometric
  spreading} along the ray tube. We denote it $J(\sigma)$. Letting
$\mathsf{A}$ denote a polynomial approximation of $A$ off of the grid
nodes in $\Omega_h$, using~\cref{eq:amplitude-update}, we can compute
$A(\xhat)$ from:
\begin{equation}
  A(\xhat) = \sqrt{\frac{c_0(L_{\lam})}{|\det(\m{Q}(L_{\lam}))|}} \mathsf{A}(\xlam) = \sqrt{\frac{c(\xhat)}{J(\xlam)}} \mathsf{A}(\xlam).
\end{equation}
Since this depends on $\m{Q}(L)$, we must
solve~\eqref{eq:hamiltons-equations-ray-tube} along $\mphi$, requiring
us to provide initial conditions at $\sigma = 0$. Note that if we set
$\sigma = 0$ in~\eqref{eq:amplitude-update}, we can see that
$|\det(\m{Q}(0))| = c_0(0)$ is necessary. A simple choice for the
initial conditions for $\m{Q}$ is $\m{Q}(0) = c_0(0)^{1/n}
\m{I}$. This assumes that we aren't too close to a point source, where
$A$ is singular.

To find initial conditions for $\m{P}$, first expand $\tau$ in a
Taylor series orthogonal to the central ray, i.e.\ in the coordinates
$\m{q}$. Doing this, we find that:
\begin{equation}
  \tau(\phi(\sigma, \m{q})) = \tau(\phi_0(\sigma)) + \frac{1}{2} \m{q}^\top \left.\frac{\partial^2 \tau}{\partial \m{q} \partial \m{q}^\top}\right|_{\m{q}=0} \m{q} + O(\m{q}^3).
\end{equation}
In this Taylor expansion, the linear term disappears since the rays
and wavefronts are orthogonal. If we let:
\begin{equation}
  \m{\Gamma} = \left.\frac{\partial^2 \tau}{\partial\m{q}\partial\m{q}^\top}\right|_{\m{q}=0},
\end{equation}
we find that $\m{\Gamma}$ satisfies the matrix Ricatti equation:
\begin{equation}\label{eq:Gamma-ricatti}
  \frac{d\m{\Gamma}}{d\sigma} + c_0 \m{\Gamma}^2 + \frac{1}{c_0^2} \left.\frac{\partial^2 c}{\partial\m{q}\partial\m{q}^\top}\right|_{\m{q}=0} = 0.
\end{equation}
The standard way to solve~\eqref{eq:Gamma-ricatti} is to use the
ansatz $\m{\Gamma} = \m{P}\m{Q}^{-1}$, which, indeed, leads us back
to~\eqref{eq:hamiltons-equations-ray-tube}. However, this viewpoint
furnishes us with the initial conditions for $\m{P}$, since
$\m{\Gamma}(0)$ can now be readily computed from
$\nabla^2\mathsf{T}(\xlam)$.

\paragraph{Marching the amplitude of a linear speed of sound} As a
simple but important test case, we consider a problem with a constant
speed of sound, i.e.:
\begin{equation}
  s(\m{x}) = \frac{1}{c(\m{x})}, \qquad c(\m{x}) = v_0 + \m{v}^\top \m{x}.
\end{equation}
In this case, \eqref{eq:hamiltons-equations-ray-tube} simplifies
considerably since $\nabla^2 c \equiv 0$, implying
$\m{P}(\sigma) \equiv \m{P}(0) = \m{\Gamma}(0) \m{Q}(0)$. From this, we can
integrate $d\m{Q}/d\sigma$ from 0 to $L$ to obtain:
\begin{equation}
  \m{Q}(L) = \Bigg[\m{I} + \bigg(\int_0^L c(\mphi(\sigma))d\sigma\bigg) \m{\Gamma}(0) \Bigg]\m{Q}(0).
\end{equation}
Denote the integral in this expression for $\m{Q}(L)$ by
$\epsilon$. To evaluate $\epsilon$ approximately, we can apply the
trapezoid rule to get:
\begin{equation}
  \epsilon = \int_0^L c(\mphi(\sigma))d\sigma = L \cdot \left(v_0 + \m{v}^\top (\xhat + \xlam)/2\right) + O(L^2),
\end{equation}
which implies that $|\epsilon| = O(L)$. The fact that the error is
$O(L^2)$ in this case follows from usual error bound for the trapezoid
rule and the fact that
$\max_{0\leq\sigma\leq L}|\mphi''(\sigma)| = O(L^{-1})$, by our choice
of parametrization.

We would like to develop a simple update rule for the geometric
spreading. First, note that the determinant satisfies the following
identity:
\begin{equation}
  \det(\m{I} + \epsilon \m{\Gamma}(0)) = 1 + \epsilon \operatorname{tr}(\m{\Gamma}(0)) + O(\epsilon^2).
\end{equation}
Next, recall that $\m{E}(0)$ is an orthogonal matrix such that
$\m{E}(0)^\top \mphi_0' = \m{E}(0)^\top \tlam = 0$. Let
$\m{U} = \begin{bmatrix} \m{E}(0) & \tlam \end{bmatrix}$ and write:
\begin{equation}\label{eq:trace-hessian}
  \operatorname{tr} \nabla^2 T(\xlam) = \operatorname{tr} \m{U}^\top \nabla^2 T(\xlam) \m{U} = \operatorname{tr} \m{E}(0)^\top \nabla^2 T(\xlam) \m{E}(0) + \tlam^\top \nabla^2 T(\xlam) \tlam.
\end{equation}
By definition,
$\m{\Gamma}(0) = \m{E}(0)^\top \nabla^2 T(\xlam) \m{E}(0)$. Taking the
gradient of~\eqref{eq:eikonal-equation}, we get:
\begin{equation}
  \nabla^2 T(\xlam) \nabla T(\xlam) = s(\xlam) \nabla s(\xlam),
\end{equation}
which leads immediately to:
\begin{equation}\label{eq:tlam-2nd-der}
  \tlam^\top \nabla^2 T(\xlam) \tlam = \tlam^\top \nabla s(\xlam),
\end{equation}
noting that $\tlam = \nabla T(\xlam)/\|\nabla
T(\xlam)\|$. Combining~\eqref{eq:trace-hessian}
and~\eqref{eq:tlam-2nd-der} gives:
\begin{equation}
  \operatorname{tr}\m{\Gamma}(0) = \Delta T(\xlam) - \tlam^\top \nabla s(\xlam),
\end{equation}
since $\operatorname{tr}\nabla^2 T(\xlam) = \Delta T(\xlam)$. This
gives the following update for $J$:
\begin{equation}\label{eq:J-update}
  J(\xhat) = \Big|1 + \epsilon \cdot \big(\Delta T(\xlam) - \tlam^\top \nabla s(\xlam)\big)\Big| \cdot \mathsf{J}(\xlam).
\end{equation}
Here, $\mathsf{J}$ denotes a local polynomial approximation to
$J$. This can be computed directly from data immediately available
after solving the optimization problem that determines $T(\xhat)$ and
$\nabla T(\xhat)$.

\paragraph{Initial data for $J$ and $A$} Determing the initial data
for the amplitude is
involved~\cite{Avila:1963aa,Babich:1979aa,Popov:2002aa,Qian:2016aa},
and detailed consideration of this problem is outside the scope of
this work. Instead, we note that for a point source in 2D, the
following hold approximately near the point source:
\begin{equation}\label{eq:J-and-A-initial-data}
  J(\m{x}) \sim |\m{x}|, \qquad A(\m{x}) \sim \frac{e^{i \pi/4}}{2 \sqrt{2\pi\omega}} \sqrt{\frac{c(\m{x})}{J(\m{x})}}.
\end{equation}
See Popov for quick derivations of these
approximations~\cite{Popov:2002aa}. In our test problems, we
initialize $J$ to $|\m{x}|$ near the point source, march $J$ according
to~\eqref{eq:J-update} where
$\mathsf{J} = (1 - \lambda) J(\m{x}_1) + \lambda J(\m{x}_2)$, and
compute the final amplitude from:
\begin{equation}\label{eq:amplitude-formula}
  A(\m{x}) = \frac{e^{i \pi/4}}{2\sqrt{2\pi\omega}} \sqrt{\frac{c(\m{x})}{J(\m{x})}}.
\end{equation}
We emphasize that this is only valid for two-dimensional problems. The
same sort of approach can be used for 3D problems,
but~\eqref{eq:J-and-A-initial-data} must be modified.

\paragraph{Marching the amplitude for more general slowness functions}
The update given by~\eqref{eq:J-update} is valid if we approximate the
speed function $c = 1/s$ with a piecewise linear function with nodal
values taken from $c(\m{x})$, where $\m{x} \in \Omega_h$. This should
be a reasonable thing to do, since the update rule given
by~\eqref{eq:J-update} in this case appears to be $O(h^2)$
accurate. Since the accuracy of $\nabla^2 T$ computed by our method is
limited, we should not expect to be able to obtain much better than
$O(h)$ accuracy for $J$. That said, a more accurate update for $J$
could be obtained by numerically
integrating~\eqref{eq:hamiltons-equations-ray-tube}.

\section{Numerical experiments}\label{sec:numerical-results} In this section, we first present a
variety of test problems which differ primarily in the choice of
slowness function $s$. The choices of $s$ range from simple, such as
$s \equiv 1$ (an overly simplified but reasonable choice for speed of
sound in room acoustics), to more strongly varying. We then present
experimental results for our different JMMs as applied to these
different slowness functions, demonstrating the significant effect the
choice of $s$ has on solver accuracy. The solvers used in these
experiments are:
\begin{itemize}
\item \texttt{JMM1}: $\mphi$ is approximated using a cubic curve, and
  tangent vectors are found by
  solving~\ref{eq:minimization-problem-1}.
\item \texttt{JMM2}: $\mphi$ is approximated using a cubic curve, with
  $\that$ optimized from~\ref{eq:minimization-problem-2} and $\tlam$
  found from Hermite interpolation at the base of the update.
\item \texttt{JMM3}: $\mphi$ is approximated using a quadratic curve,
  with its tangent vectors being found by optimizing.
\item \texttt{JMM4}: \texttt{JMM2} combined with the cell-marching
  method described in section~\ref{sec:cell-marching}.
\end{itemize}
We also plot the same results obtained by the
FMM~\cite{Sethian:1996aa} and
\texttt{olim8\_mp0}~\cite{Potter:2019ab}. We do not include least
squares fits for these solvers in our tables. They are mostly $O(h)$,
with some exceptions for $\nabla T$ as computed by the FMM.

We note that $s$ does not significantly affect the
runtime of any of our solvers---formally, our solvers run in
$O(|\Omega_h| \log |\Omega_h|)$ time, where the constant factors are
essentially insensitive to the choice of $s$. We note that the cost of
updating the heap is very small compared to the cost of doing
updates. Since only $|\Omega_h|$ updates must be computed, the CPU
time of the solver effectively scales like $O(|\Omega_h|)$ for all
problem sizes considered in this paper.

\subsection{Test problems} In this section, we provide details for the
test problems used in our numerical tests.

\paragraph{Constant slowness with a point source} For this problem, the slowness and solution are given by:
\begin{equation}
  s \equiv 1, \qquad \tau(\m{x}) = \|\m{x}\|.
\end{equation}
We take the domain to be
$\Omega = [-1, 1] \times [-1, 1] \subseteq \mathbb{R}^2$. To control
the size of the discretized domain, we let $M > 0$ be an integer and
set $h = 1/M$, from which we define $\Omega_h$ accordingly. We place a
point source at $\m{x}^\circ = (0, 0) \in \Omega_h$. The set of
initial boundary data locations given by is
$\Gamma_h = \{\m{x}^\circ\}$, with boundary conditions given by
$g(\m{x}^\circ) = 0$.

\paragraph{Linear speed with a point source (\#1)} Our next test
problem has a linear velocity profile. This might model the variation
in the speed of sound due to a linear temperature gradient (e.g., in a
large room). The slowness is given
by~\cite{Fomel:2009aa,slotnick1959lessons}:
\begin{equation}
  s(\m{x}) = \left[\frac{1}{s_0} + \m{v}^\top \m{x}\right]^{-1},
\end{equation}
where $s_0 > 0$, and $\m{v} \in \mathbb{R}^2$ are parameters. The
solution is given by:
\begin{equation}\label{eq:eikonal-linear-speed}
  \tau(\m{x}) = \frac{1}{\|\m{v}\|} \cosh^{-1}\hspace{-0.1em}\left(1 + \frac{1}{2} s_0 s(\m{x}) \|\m{v}\|^2 \|\m{x}\|^2\right).
\end{equation}
For our first test with a linear speed function, we take $s_0 = 1$ and
$\m{v} = (0.133, -0.0933)$. For this problem,
$\Omega = [-1, 1] \times [-1, 1]$, $\Gamma_h = \{\m{x}^\circ\}$, and $g(\m{x}^\circ) = 0$.

\paragraph{Linear speed with a point source (\#2)} For our second
linear speed test problem, we set $s_0 = 2$ and $\m{v} = (0.5, 0)$ as
in~\cite{Qi:2019aa}. For this problem, we let
$\Omega = [0, 1] \times [0, 1]$, discretize into $M$ nodes along each
axis, and define $\Omega_h$ accordingly (i.e., $|\Omega_h| = M^2$,
with $M = h^{-1}$). We take $\m{x}^\circ$, $\Gamma_h$, and $g$ to be
same as in the previous two test problems.

\paragraph{A nonlinear slowness function involving a sine function}
For $\m{x} = (x_1, x_2)$, we set:
\begin{equation}
  \tau(\m{x}) = \frac{x_1^2}{2} + 2 \sin\hspace{-0.1em}\left(\frac{x_1 + x_2}{2}\right)^2.
\end{equation}
This eikonal has a unique minimum, $\tau(0, 0) = 0$, and is strictly
convex in $\Omega = [-1, 1] \times [-1, 1]$. This lets us determine the slowness
from the eikonal equation, giving:
\begin{equation}
  s(\m{x}) = \sqrt{\sin(x_1 + x_2)^2 + \big(x_1 + \sin(x_1 + x_2)\big)^2}.
\end{equation}
For this test problem, we take $\Gamma_h$ and $\Omega_h$ as in the
constant slowness point source problem.

\paragraph{Sloth} A slowness function called ``sloth'' (jargon from
geophysics) is taken from Example 1 of Fomel et
al.~\cite{Fomel:2009aa}:
\begin{equation}
  s(\m{x}) = \sqrt{s_0^2 + 2 \m{v}^\top \m{x}}.
\end{equation}
For our test with this slowness function, we set $s_0 = 2$, and
$\m{v} = (0, -3)$. In this case, to avoid shadow zones formed by caustics, we take
$\Omega = [0, \tfrac{1}{2}] \times [0, \tfrac{1}{2}]$. The discretized domain and boundary
data are determined analogously to the earlier cases.

\begin{figure}
  \centering
  \includegraphics[width=\linewidth]{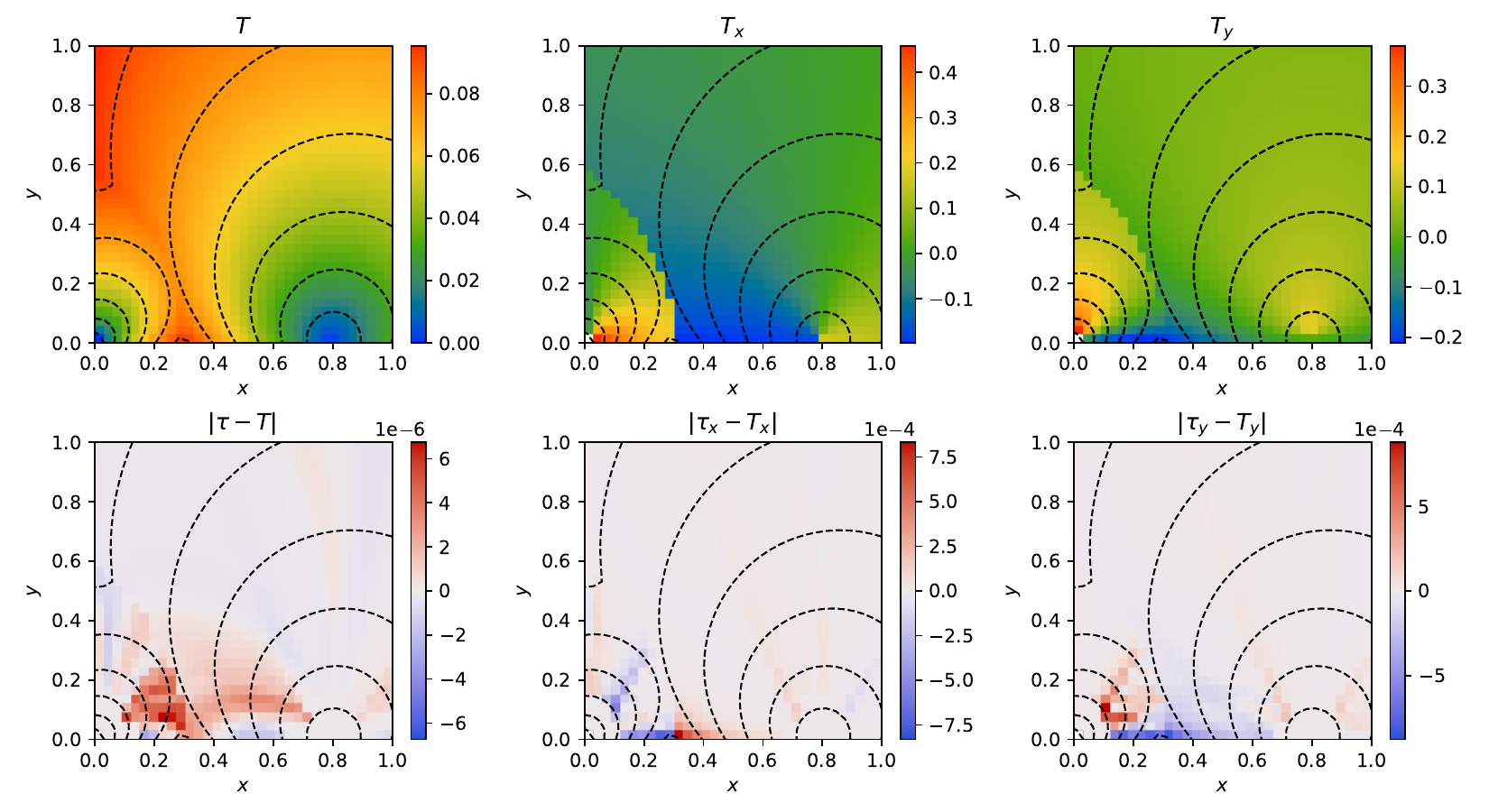}
  \vspace{-2em}
  \caption{A test problem with two point sources. The domain is
    $\Omega = [0,1]^2$ discretized into a $33 \times 33$ grid. The top
    row contains computed values and the bottom row contains signed
    errors. We overlay several contours of the numerically computed
    eikonal, $T$, for context. The two point sources are separated by
    a shockline, which is accurately localized even for this extremely
    coarse mesh. E.g., if we compared values of $T_x$ or $T_y$ on the
    ``wrong'' side of the shockline, we would observe $O(1)$
    error. \emph{Left}: $T$. \emph{Middle}: $T_x$. \emph{Right}:
    $T_y$. }\label{fig:two-point-sources}
\end{figure}

\paragraph{Two point sources} We consider a linear speed function:
\begin{equation}
  c(\m{x}) = \frac{1}{s(\m{x})} = 2 + 5x_1 + 20x_2,
\end{equation}
with point sources located at $\m{x}_0 = (0, 0)$ and
$\m{x}_1 = (0.8, 0)$. This is the same problem considered in Figure 4
of Qi and Vladimirsky~\cite{Qi:2019aa}. We use
\eqref{eq:eikonal-linear-speed} to compute the groundtruth eikonal for
each point source. If we let $\tau_0$ and $\tau_1$ denote the eikonals
for each point source problem considered individually, then the
combined eikonal is:
\begin{equation}
  \tau(\m{x}) = \min(\tau_0(\m{x}), \tau_1(\m{x})).
\end{equation}
We can easily see that all derivatives of $\tau$ are undefined on the
shockline, where $\tau_0 \equiv \tau_1$. However, away from the
shockline, the derivatives are well-defined and smooth: if
$j = \Arg\min_i\tau_i(\m{x})$, then $D\tau \equiv D\tau_j$. We solve
this problem on $\Omega = [0, 1]^2$, so that $\Omega_h$ is a grid with
$N = 33$ nodes in each direction. The results are shown in
\Cref{fig:two-point-sources}.

For this problem, we note that the shockline is localized sufficiently
well: all nodes $\m{x} \in \Omega_h$ are on the correct side, so that
when we compute the errors, each nodal value is compared with the
correct branch of $\tau = \min(\tau_0, \tau_1)$. Note the use of an
extremely coarse mesh. The mesh used here is coarser than any mesh
used in Qi and Vladimirsky's Figure 4, but the eikonal achieves a
smaller maximum error than all of their test problems.

\begin{figure}
  \centering
  \includegraphics[width=\linewidth]{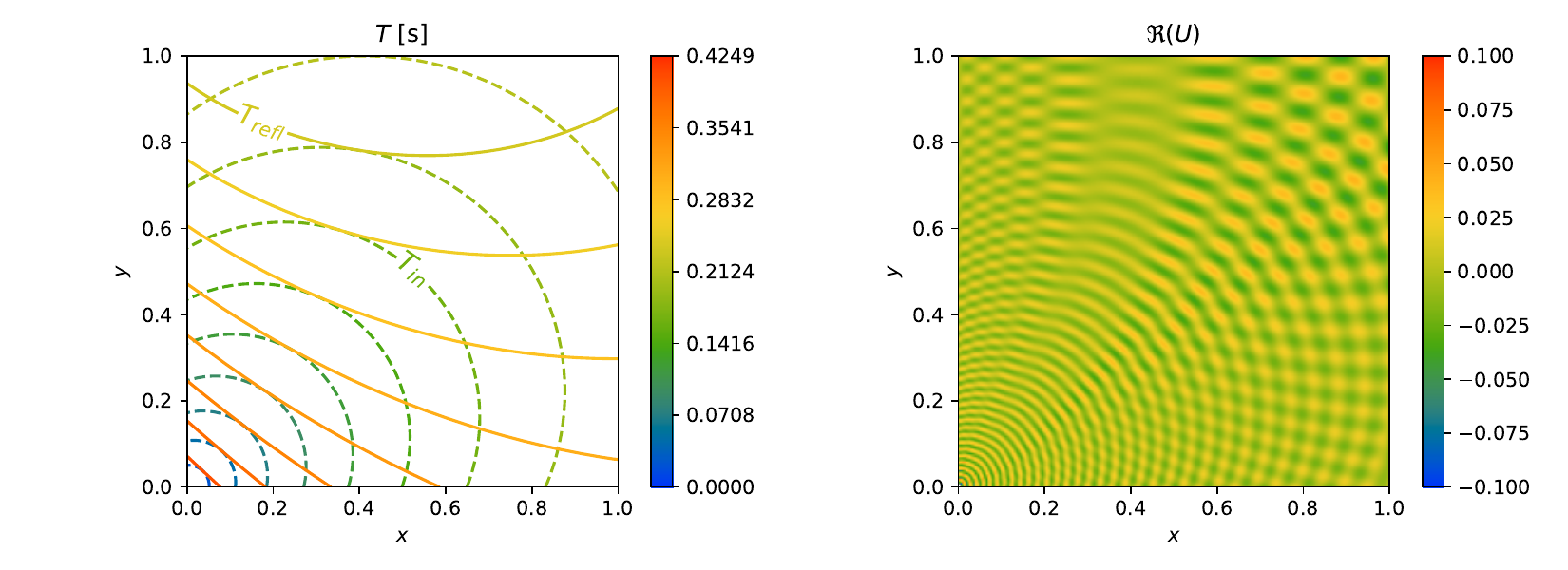}
  \vspace{-2em}
  \caption{A single reflection from a wall in a semi-infinite domain
    for $\omega = 1000$ (we consider only the top edge of the boundary
    to be reflecting). \emph{Left}: $T$ for the incident and reflected
    fields. The reflected field satisfies a specular reflection
    condition along the edge of the domain. \emph{Right}: the real
    part of the approximation to the solution of the Helmholtz
    equation so obtained, given by
    \eqref{eq:reflection-equation-for-example}.}\label{fig:reflection}
\end{figure}

\paragraph{A single reflection} We additionally include a simple test
for computing multiple arrivals. For the linear speed function:
\begin{equation}
  c(\m{x}) = \frac{1}{s(\m{x})} = 2 + 5x_1 + 7x_2,
\end{equation}
we solve \eqref{eq:eikonal-equation} on $\Omega = [0, 1]^2$,
discretized into $N = 101$ nodes in each direction. We place a point
source at $\m{x} = (0, 0)$ and compute the eikonal, which we denote
$\tau_{\operatorname{in}}$. We then restrict
$\tau_{\operatorname{in}}$ and $\nabla \tau_{\operatorname{in}}$
(after reflection) to the set $\Gamma = [0, 1] \times \{1\}$ (the top
edge of the domain), and solve the reflected eikonal equation:
\begin{align}
  \|\nabla \tau_{\operatorname{refl}}(\m{x})\| &= s(\m{x}), &\m{x} \in \Omega, \\
  \tau_{\operatorname{refl}}(\m{x}) &= \tau_{\operatorname{in}}(\m{x}), &\m{x} \in \Gamma, \\
  \left.\frac{\partial\tau_{\operatorname{refl}}}{\partial x}\right|_{\m{x}} &= \left.\frac{\partial\tau_{\operatorname{in}}}{\partial x}\right|_{\m{x}}, &\m{x} \in \Gamma, \\
  \left.\frac{\partial\tau_{\operatorname{refl}}}{\partial y}\right|_{\m{x}} &= -\left.\frac{\partial\tau_{\operatorname{in}}}{\partial y}\right|_{\m{x}}, &\m{x} \in \Gamma.\label{eq:reflected-tau-y-BC}
\end{align}
Note the minus sign in \eqref{eq:reflected-tau-y-BC}. This corresponds
to a specular reflection from the ``wall'' $\Gamma$. After we compute
$T_{\operatorname{in}}$ and $T_{\operatorname{refl}}$ numerically, we
can then compute the geometric spreading $J_{\operatorname{in}}$ and
reflected geometric spreading $J_{\operatorname{refl}}$. Since the
reflecting set is flat, we can use the boundary
condition~\cite{Popov:2002aa}:
\begin{equation}
  J_{\operatorname{refl}}(\m{x}) = J_{\operatorname{in}}(\m{x}), \qquad \m{x} \in \Gamma.
\end{equation}
Afterwards, we use \eqref{eq:amplitude-formula} to obtain:
\begin{equation}\label{eq:reflection-equation-for-example}
  U(\m{x}) = A_{\operatorname{in}}(\m{x})\exp(-i\omega T_{\operatorname{in}}(\m{x})) + A_{\operatorname{refl}}(\m{x})\exp(-i\omega T_{\operatorname{refl}}(\m{x})).
\end{equation}
In \Cref{fig:reflection}, we plot $T_{\operatorname{in}}$,
$T_{\operatorname{out}}$, and the real part of $U$.

\subsection{Experimental results}

\begin{figure}
  \includegraphics[width=\textwidth]{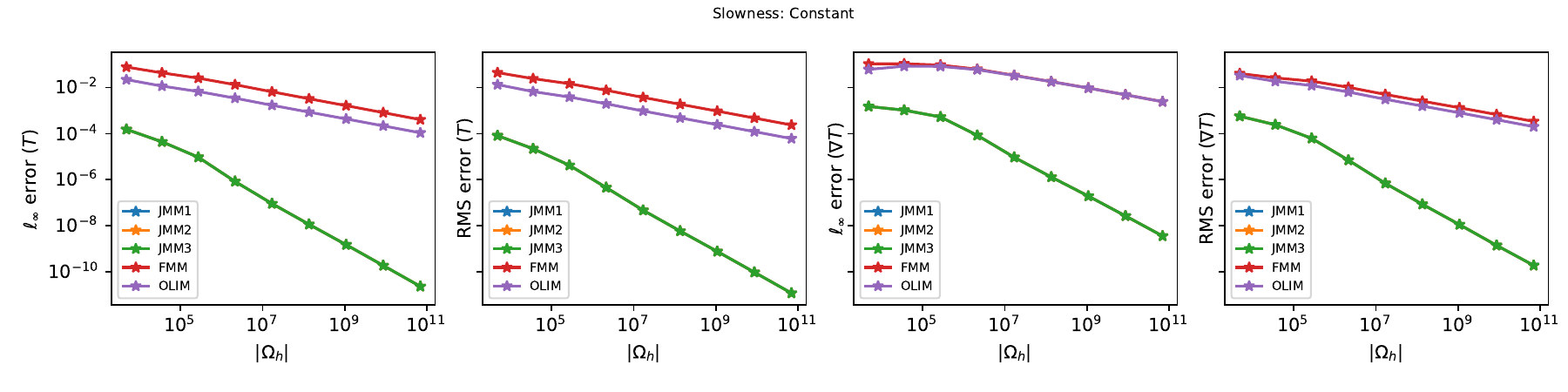}
  \includegraphics[width=\textwidth]{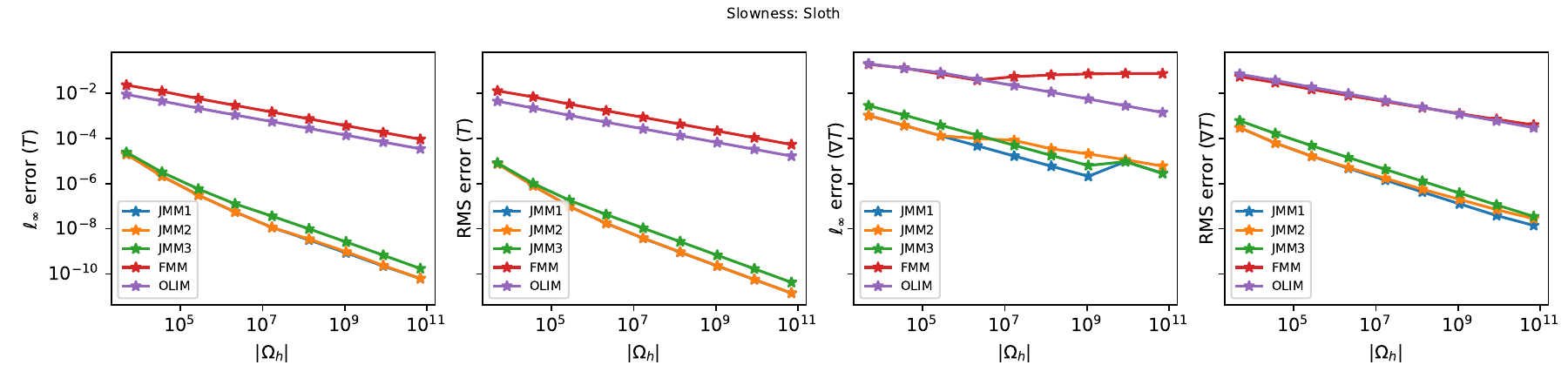}
  \includegraphics[width=\textwidth]{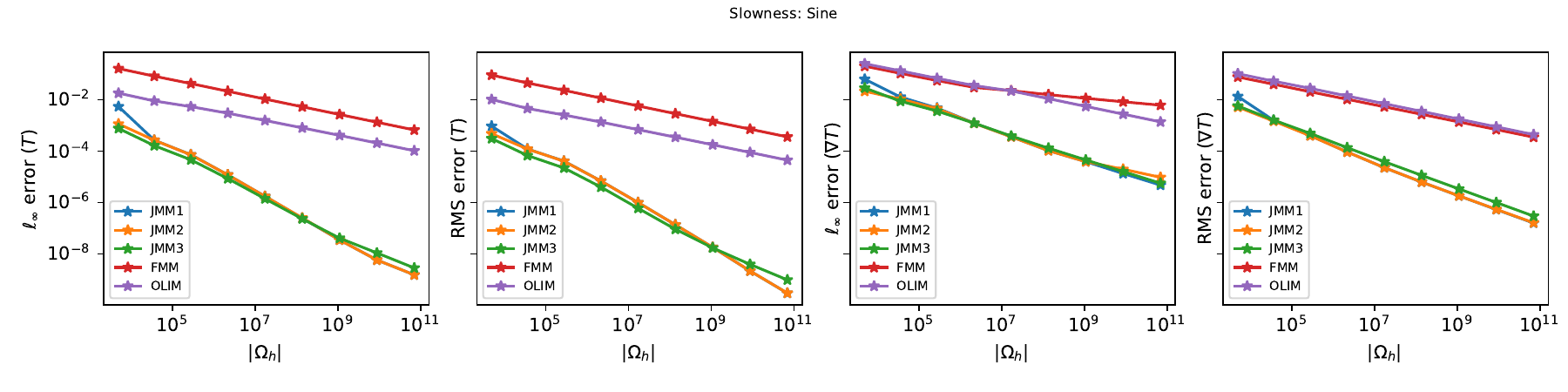}
  \includegraphics[width=\textwidth]{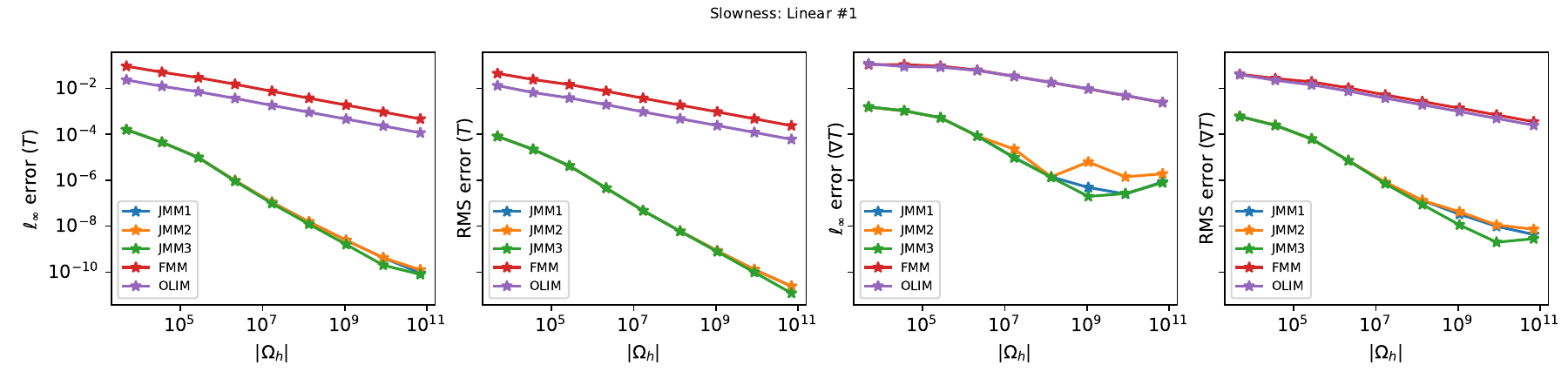}
  \includegraphics[width=\textwidth]{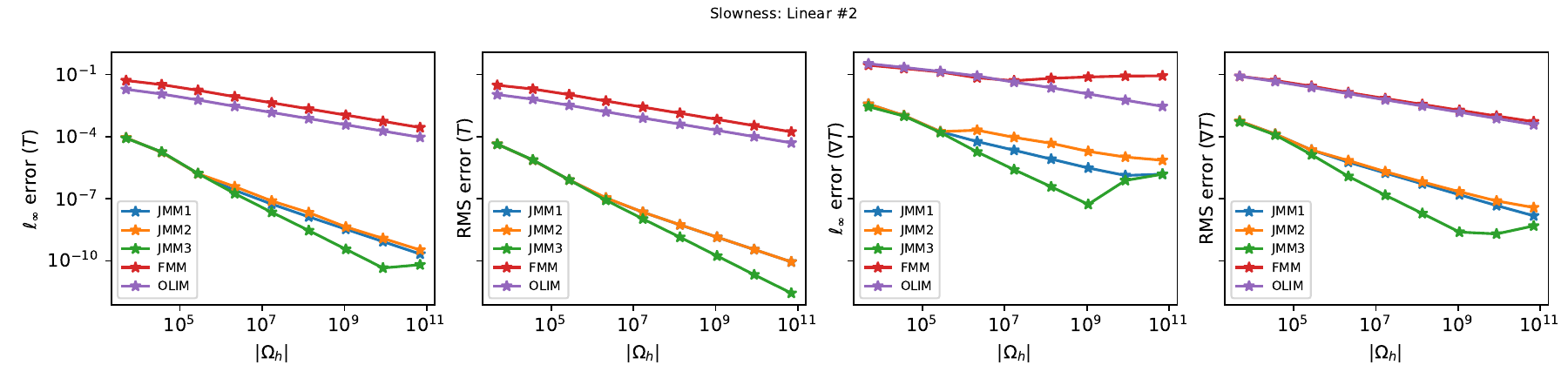}
  \vspace{-2em}
  \caption{Plots comparing domain size ($|\Omega_h|$) and
    $\ell_\infty$ and RMS errors for $T$ and
    $\nabla T$.}\label{fig:size-vs-error}
\end{figure}

\begin{figure}
  \includegraphics[width=\textwidth]{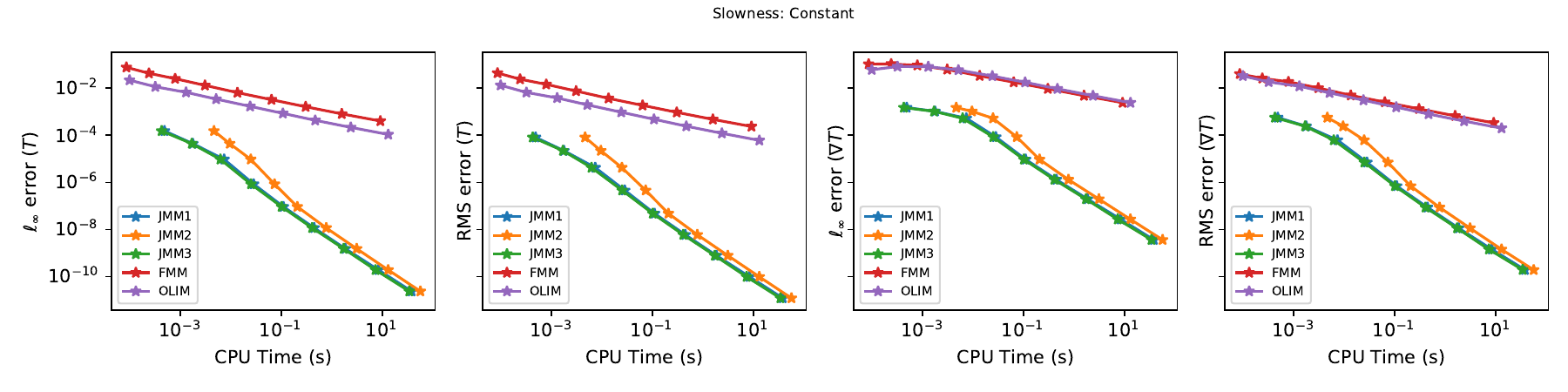}
  \includegraphics[width=\textwidth]{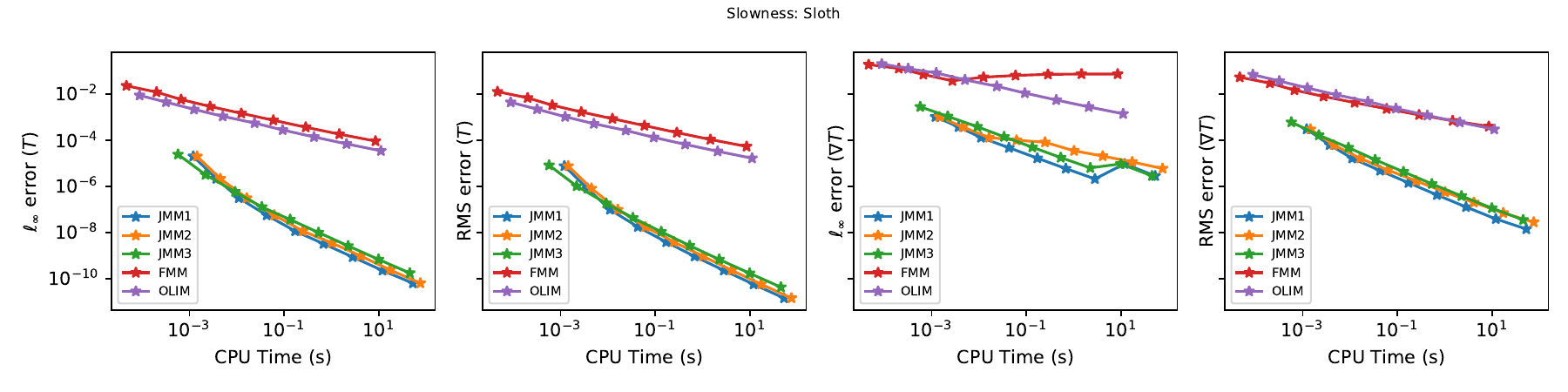}
  \includegraphics[width=\textwidth]{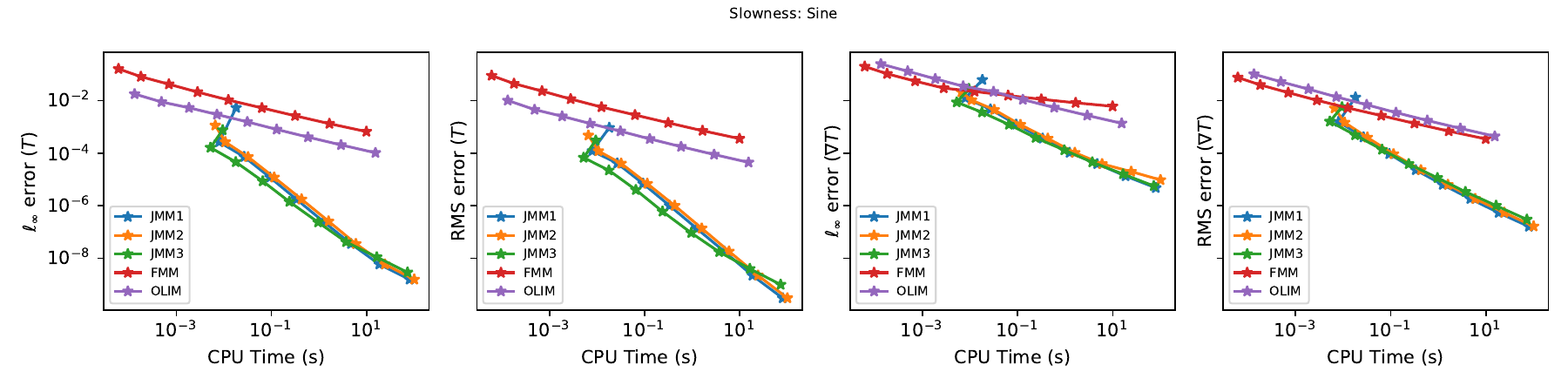}
  \includegraphics[width=\textwidth]{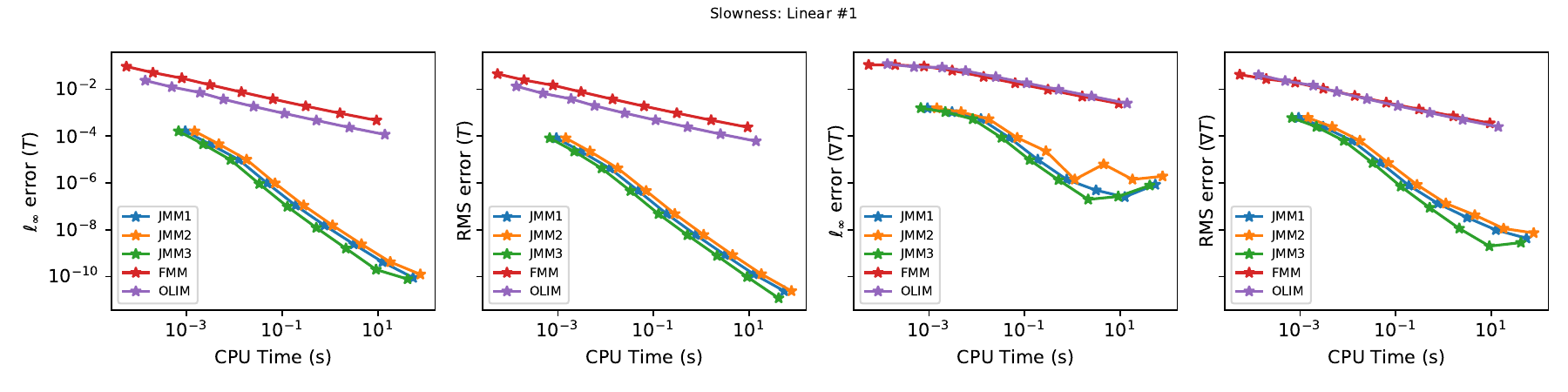}
  \includegraphics[width=\textwidth]{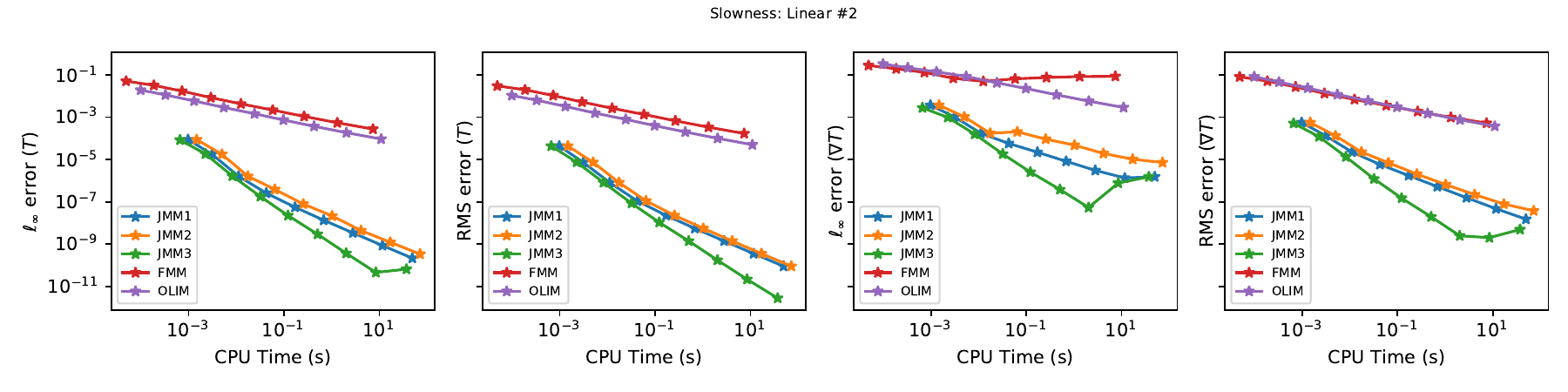}
  \vspace{-2em}
  \caption{Plots comparing CPU runtime in seconds and
    errors.}\label{fig:time-vs-error}
\end{figure}

\begin{figure}
  \includegraphics[width=\textwidth]{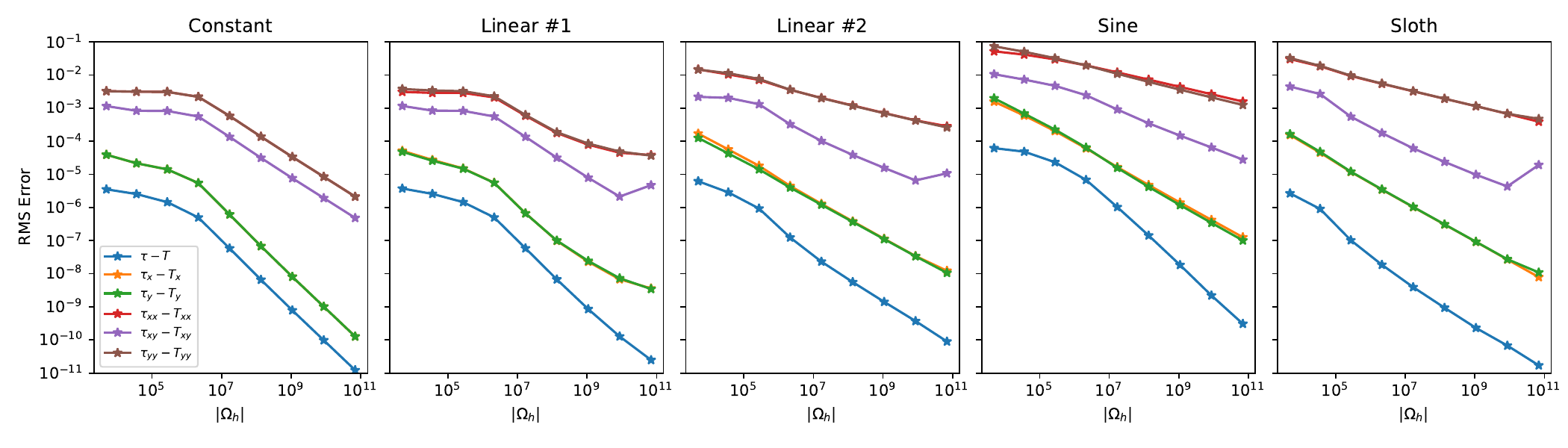}
  \vspace{-2em}
  \caption{Domain size vs.\ RMS error for
    \texttt{JMM4}.}\label{fig:size-vs-error-jmm4}
\end{figure}

\begin{table}
  \centering
  \begin{tabular}{cccccc}
    & JMM & $E_{\mbox{max}} (T)$ & $E_{\mbox{RMS}} (T)$& $E_{\mbox{max}} (\nabla T)$ & $E_{\mbox{RMS}} (\nabla T)$ \\
    \midrule
    \multirow{3}{*}{Constant}
    & \#1 & 2.87 & 2.87 & 2.28 & 2.72 \\
    & \#2 & 2.87 & 2.87 & 2.28 & 2.72 \\
    & \#3 & 2.87 & 2.87 & 2.28 & 2.72 \\
    \midrule
    \multirow{3}{*}{Linear \#1}
    & \#1 & 2.77 & 2.85 & 2.14 & 2.52 \\
    & \#2 & 2.77 & 2.85 & 1.70 & 2.48 \\
    & \#3 & 2.86 & 2.87 & 2.28 & 2.73 \\
    \midrule
    \multirow{3}{*}{Linear \#2}
    & \#1 & 2.48 & 2.52 & 1.70 & 1.97 \\
    & \#2 & 2.38 & 2.52 & 1.16 & 1.88 \\
    & \#3 & 3.03 & 3.03 & 2.70 & 3.02 \\
    \midrule
    \multirow{3}{*}{Sine}
    & \#1 & 2.76 & 2.57 & 1.77 & 2.09 \\
    & \#2 & 2.51 & 2.46 & 1.58 & 1.94 \\
    & \#3 & 2.37 & 2.38 & 1.54 & 1.79 \\
    \midrule
    \multirow{3}{*}{Sloth}
    & \#1 & 2.39 & 2.48 & 1.49 & 1.84 \\
    & \#2 & 2.37 & 2.47 & 0.87 & 1.73 \\
    & \#3 & 2.15 & 2.21 & 1.47 & 1.76 \\
  \end{tabular}
  \vspace{1em}
  \caption{The order of convergence $p$ for each combination of test
    problems and solvers, computed for different types of errors and
    fit as $C h^p$.}\label{fig:least-squares-fits}
\end{table}

\begin{table}
  \centering
  \begin{tabular}{ccccccc}
    & $\tau - T$ & $\tau_x - T_x$ & $\tau_y - T_y$ & $\tau_{xx} - T_{xx}$ & $\tau_{xy} - T_{xy}$ & $\tau_{yy} - T_{yy}$ \\
    \midrule
    Constant & 3.09 & 3.11 & 3.11 & 2.01 & 2.05 & 2.01 \\
    Linear \#1 & 2.99 & 2.43 & 2.40 & 1.39 & 2.01 & 1.39 \\
    Linear \#2 & 2.10 & 1.76 & 1.72 & 0.77 & 1.25 & 0.77 \\
    Sine & 2.91 & 1.80 & 1.89 & 0.73 & 1.31 & 0.80 \\
    Sloth & 2.03 & 1.76 & 1.75 & 0.75 & 1.33 & 0.76 \\
  \end{tabular}
  \caption{The order of convergence $p$ for \texttt{JMM4} for each
    component of the total 2-jet of $\tau$, computed from least
    squares fits of the RMS error. The fits only incorporate the 4th
    through the 8th problem sizes to avoid artifacts for small and
    large problem sizes. See
    Figure~\ref{fig:size-vs-error-jmm4}.}\label{fig:least-squares-fits-jmm4}
\end{table}

\begin{figure}
  \centering
  \includegraphics[width=\linewidth]{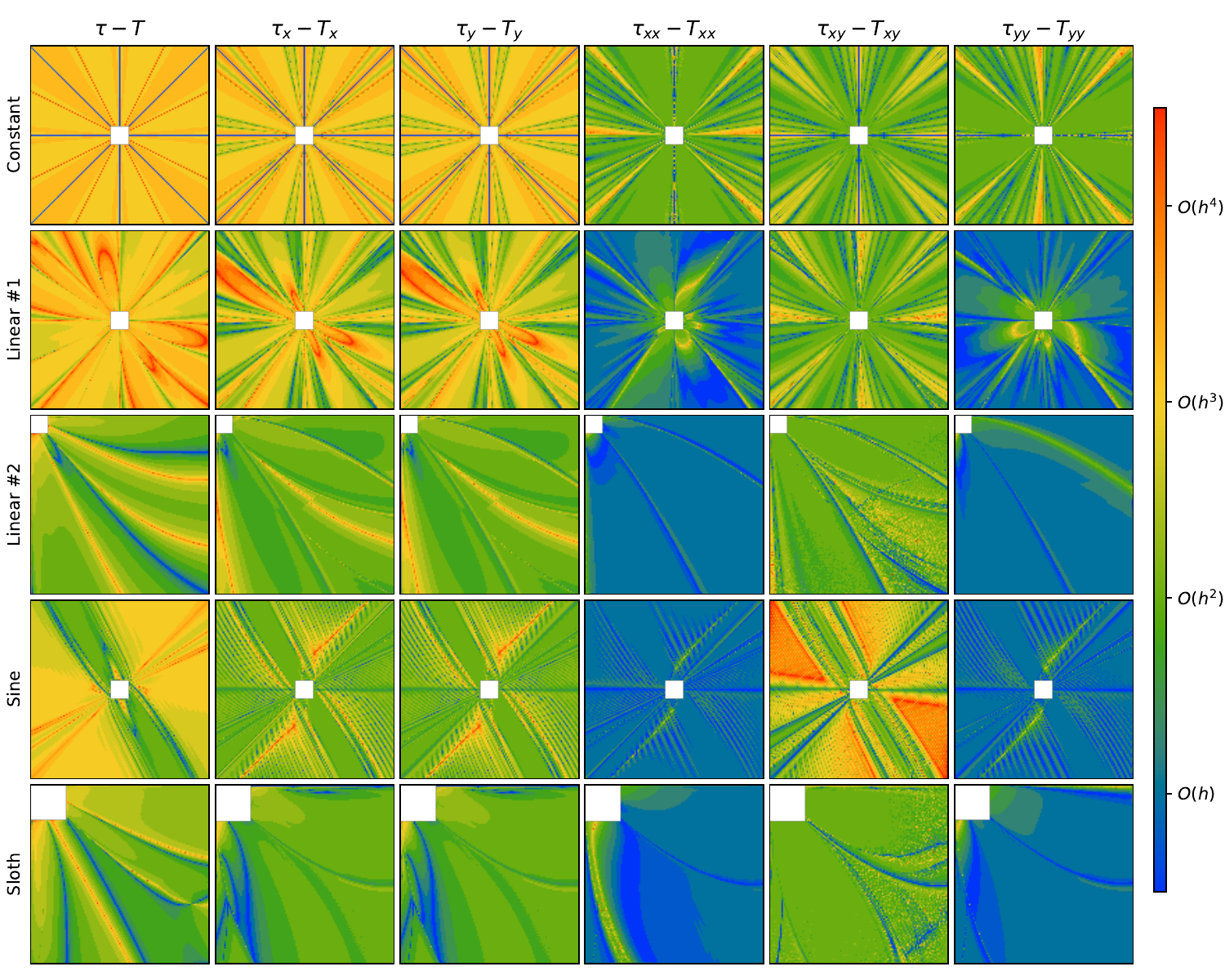}
  \vspace{-2em}
  \caption{A plot of the pointwise convergence at each point in
    $\Omega_h$ for \texttt{JMM4}. To obtain these plots, starting with
    $N = 129$, we decimate each larger problem size (up to
    $N = 2,049$) to a $129 \times 129$ grid, and do a least squares
    fit at each point. This gives us an estimate of the order of
    convergence at each point.}\label{fig:jmm4-pointwise-convergence}
\end{figure}

\begin{figure}
  \centering
  \includegraphics[width=\linewidth]{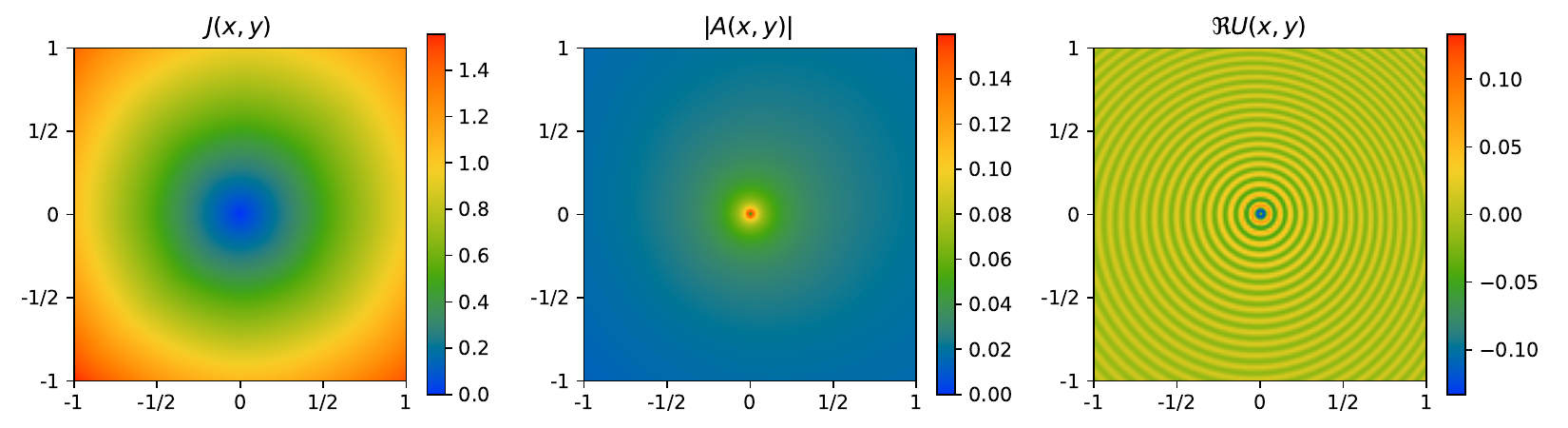}
  \vspace{-2em}
  \caption{Plots related to computing the amplitude and a numerical
    approximation to the solution to
    $(\Delta + \omega^2 s(\m{x})^2) u(\m{x}) = \delta(\m{x})$, denoted
    $U(\m{x})$ for the Linear \#1 test problem. \emph{Left}: the
    geometric spreading. \emph{Middle}: the amplitude
    function. \emph{Right}: the numerical solution
    $U$.}\label{fig:p-plots}
\end{figure}

The results of our numerical experiments evaluating the JMMs described
in Section~\ref{sec:minimization-problems} are presented in
Table~\ref{fig:least-squares-fits} and
Figures~\ref{fig:size-vs-error} and~\ref{fig:time-vs-error}. The
numerical tests for \texttt{JMM4}, which uses cell marching, are given
in Table~\ref{fig:least-squares-fits-jmm4} and
Figures~\ref{fig:size-vs-error-jmm4}
and~\ref{fig:jmm4-pointwise-convergence}. An example where the
geometric spreading and amplitude are computed using cell marching
method is shown in Figure~\ref{fig:p-plots}.

For more benign choices of $s$, the errors generally convergence with
$O(h^3)$ accuracy for $T$ and $O(h^2)$ accuracy for $\nabla T$ in the
RMS error. For the special case of $s \equiv 1$, the gradients also
converge with nearly $O(h^3)$ accuracy. For more challenging nonlinear
choices of $s$, the eikonal converges with somewhere between $O(h^2)$
and $O(h^3)$ accuracy, while the gradient converges with nearly
$O(h^2)$ accuracy.

We note that in some cases the gradient begins to diverge for large
problem sizes. This occurs because our tolerance for minimizing $F$ is
not small enough, and also because $\nabla^2 F$ is $O(h)$. For our
application, our goal is to save memory and compute time by using a
higher-order solver; it is unlikely we would solve problems with such
a fine discretization in practice. At the same time, choosing the
tolerance for numerical minimization based on $h$ is of
interest---partly to see how much time can be saved for coarser
problems, but also to determine to what extent the full order of
convergence can be maintained using different floating point
precisions.

The JMMs using cubic approximations for $\mphi$ tend to perform better
than those using quadratic approximations when $s$ is nonlinear and
the characteristics of~\eqref{eq:eikonal-equation} are not circular
arcs. When $s$ corresponds to a linear speed of sound, the JMMs with
quadratic $\mphi$ are a suitable choice, generally outperforming the
``cubic $\mphi$'' solvers, exhibiting cubically (or nearly cubically)
convergent RMS errors in $T$ and $\nabla T$. This is a useful finding
since the simplified solver requires fewer floating-point operations
per update, and since linear speed of sound profiles (e.g., as a
function of a linear temperature profile) are a frequently occurring
phenomenon in room acoustics.

\section{Online Package}

To recreate our results, to experiment with these solvers, and to
understand their workings, a package has been made available online on
GitHub at
\texttt{https://github.com/sampotter/jmm/tree/jmm-sisc-figures}. Details
explaining how to obtain this package and the collect the results are
available at this link.

\section{Conclusion}

We have presented a family of semi-Lagrangian label-setting methods
(\`{a} la the fast marching method) which are high-order and compact,
which we refer to as jet marching methods (JMMs). We examine a variety
of approaches to formulating one of these solvers, and in 2D, provide
extensive numerical results demonstrating the efficacy of these
approaches. We show how a form of ``adaptive'' cell-marching can be
done which is compatible with our stencil compactness requirements,
although this scheme no longer displays optimal locality.

Our solvers are motivated by problems involving repeatedly solving the
eikonal equation in complicated domains where:
\begin{itemize}
\item time and memory savings via the use of high-order solvers,
\item high-order local knowledge of characteristic directions,
\item and compactness of the solver's ``stencil'' (the neighborhood
  over which the semi-Lagrangian updates require information)
\end{itemize}
is paramount. In particular, our goal is to parametrize the multipath
eikonal in a complicated polyhedral domain in a work-efficient
manner. This solver is a necessary ingredient for carrying out this
task.

We will be continue to work along the following directions:
\begin{itemize}
\item Extension to regular grids in 3D, which should be
  straightforward and yield considerable savings over existing
  approaches, and extension to unstructured simplex meshes in 2D and
  3D. Especially in 3D, this problem is more complicated, requiring
  the computation of ``causal
  stencils''~\cite{Kimmel:1998aa,Sethian:2000aa}.
\item A rigorous proof of convergence for the solvers developed in
  this work, including a careful investigation of the conditions
  resulting in cubic convergence for both the eikonal and its gradient
  as observed in the case of constant and linear speed functions.
\end{itemize}

\section{Acknowledgments}

This work was partially supported by NSF Career Grant DMS1554907 and
MTECH Grant No.\ 6205. We thank Prof.\ Ramani Duraiswami for the
illuminating discussions and trenchant observations provided
throughout the course of this work.

\bibliographystyle{siamplain}
\bibliography{jmm}

\end{document}